  \documentclass[12pt]{amsart} 
  \usepackage{latexsym} 
  \usepackage[all]{xy} 
  \usepackage{amsfonts} 
  \usepackage{amsthm} 
  \usepackage{amsmath} 
  \usepackage{amssymb} 
  \newcommand{\cC}{{\mathcal C}} 

  \def\sw#1{{\sb{(#1)}}} 
  \def\su#1{{\sp{[#1]}}} 
   \def\suc#1{{\sp{(#1)}}} 
\def\tens{\mathop{\otimes}} 
   
  \def\endproof{\hbox{$\sqcup$}\llap{\hbox{$\sqcap$}}\medskip} 
  \def\<{{\langle}} 
  \def\>{{\rangle}} 
   
  \def\eps{\varepsilon}

  \def\note#1{{}} 
  \def\can{{\rm can}} 
   
  \def\note#1{} 
  \def\M{{\bf M}}

  \def\cM{{\bf M}}

  \def\cC{{\mathcal C}} 
   
  \def\cD{{\mathcal D}}

  \def\lrhom#1#2#3#4{{{\rm Hom}\sb{#1, #2}(#3,#4)}} 
  \def\lhom#1#2#3{{{\rm Hom}\sb{#1-}(#2,#3)}} 
  \def\rhom#1#2#3{{{\rm Hom}\sb{-#1}(#2,#3)}} 
\def\hom#1#2#3{{{\rm Hom}\sb{#1}(#2,#3)}} 
  \def\lend#1#2{{{\rm End}\sb{#1-}(#2)}} 
  \def\rend#1#2{{{\rm End}\sb{-#1}(#2)}} 
 
  \def\Lend#1#2{{{\rm End}\sp{#1-}(#2)}} 
  \def\Rend#1#2{{{\rm End}\sp{-#1}(#2)}} 
   
  \def\Rhom#1#2#3{{{\rm Hom}\sp{-#1}(#2,#3)}} 
\def\Lhom#1#2#3{{{\rm Hom}\sp{#1-}(#2,#3)}}

  \def\C{{\bf C}}

  \def\can{{\rm can}}

  \def\beq{\begin{equation}} 
  \def\eeq{\end{equation}} 
  \def\DC{{\Delta_\cC}} 
  \def \eC{{\eps_\cC}} 
  \def\DD{{\Delta_\cD}} 
  \def \eD{{\eps_\cD}}

  \def\im{{\rm Im}}

  \def\ut{{\otimes}} 
  \def\ot{{\otimes}}

  \def\roM{\varrho^{M}} 
\def\Nro{{}^{N}\varrho} 
  \newcommand{\Lra}{\Leftrightarrow} 
  \newcommand{\Ra}{\Rightarrow} 
  \def\hphi{\widehat{\varphi}}
 \def\tcan{\widetilde{\can}}
  \def\roA{\varrho^{A}}
 \def\Aro{{}^{A}\varrho}

  \newcounter{zlist} 
  \newenvironment{zlist}{\begin{list}{(\arabic{zlist})}{ 
  \usecounter{zlist}\leftmargin2.5em\labelwidth2em\labelsep0.5em 
  \topsep0.6ex
  \parsep0.3ex plus0.2ex minus0.1ex}}{\end{list}}

  \newcounter{blist} 
  \newenvironment{blist}{\begin{list}{(\alph{blist})}{ 
  \usecounter{blist}\leftmargin2.5em\labelwidth2em\labelsep0.5em 
  \topsep0.6ex 
  \parsep0.3ex plus0.2ex minus0.1ex}}{\end{list}} 

  \newcounter{rlist}

  \def\Label#1{\label{#1}\ifmmode\llap{[#1] }\else 
  \marginpar{\smash{\hbox{\tiny [#1]}}}\fi} 
  \def\Label{\label}

  \newcounter{c} 
   
  \newcommand{\etyk}[1]{\vspace{-7.4mm}$$\begin{equation}\Label{#1} 
  \addtocounter{c}{1}} 
  \renewcommand{\]}{\ifnum \value{c}=1 $$\else \end{equation}\fi} 
  \newtheorem{proposition}{Proposition}[section] 
  \newtheorem{lemma}[proposition]{Lemma} 
   
  \newtheorem{corollary}[proposition]{Corollary} 
  \newtheorem{theorem}[proposition]{Theorem} 

  \theoremstyle{definition} 
  \newtheorem{definition}[proposition]{Definition} 
  \newtheorem{example}[proposition]{Example} 

  \theoremstyle{remark}

\begin{document} 
\baselineskip 19pt

  \title{Galois comodules} 
  \author{Tomasz Brzezi\'nski} 
  \address{ Department of Mathematics, University of Wales Swansea, 
  Singleton Park, \newline\indent  Swansea SA2 8PP, U.K.} 
  \email{T.Brzezinski@swansea.ac.uk} 
  \urladdr{http//www-maths.swan.ac.uk/staff/tb} 
  \date{April 2004} 
  \subjclass{16W30, 13B02} 
  \begin{abstract} 
  Galois comodules of a coring are studied.  
  The conditions for a simple comodule to be a Galois 
comodule are found. A special class of Galois comodules termed 
principal comodules is introduced. These are defined as Galois comodules that are projective over their comodule endomorphism rings. A complete description of principal comodules in the case a background ring is a field is found. In particular it is shown that a (finitely generated and projective) right comodule of an $A$-coring $\cC$ is principal provided a lifting of the canonical map is a split epimorphism in the category of left $\cC$-comodules. This description is then used to characterise principal extensions or non-commutative principal bundles. Specifically, it is proven that, over a field, any entwining structure consisting of an algebra $A$, a coseparable coalgebra $C$ and a bijective entwining map $\psi$ together with a group-like element  in $C$ give rise to a principal extension provided the lifted canonical map is surjective. Induction of Galois and principal comodules via morphisms of corings 
is described.  A connection between the relative injectivity of a Galois comodule and the
properties of the extension of endomorphism rings associated to this comodule is revealed.
  \end{abstract} 
  \maketitle 

  \section{Introduction} 
In an attempt to achieve a better conceptual understanding of a generalisation of a 
Hopf-Galois extension known as a coalgebra-Galois extension, the notion of a 
{\em Galois coring} has been introduced in \cite{Brz:str}, and  recently investigated in 
\cite{Cae:Gal}, \cite{Wis:Gal}. Following a classical route in the ring extension theory 
along which the properties of an extension are encoded in properties of a module 
(cf.\ \cite{Sug:not}), 
it has been realised in \cite{KaoGom:com} that the proper framework for studying Galois corings is 
provided by a certain class of its comodules, known as {\em Galois comodules}. The most important 
result about Galois comodules is the Galois comodule structure theorem formulated in  \cite[Theorem~3.2]{KaoGom:com} (see Theorem~\ref{thm.galoisc} below) that incorporates the 
Galois Coring Theorem in \cite[Theorem~5.6]{Brz:str}, from which in turn one of Schneider's 
theorems on the structure of Hopf-Galois extensions \cite[Theorem~3.7]{Sch:pri} can be deduced.

The aim of the present paper is to study properties of Galois comodules, as a means for providing 
deeper conceptual understanding of Galois-type extensions, in particular those that are motivated 
by non-commutative geometry (where they appear as non-commutative principal bundles). In 
particular, Schneider's theorem \cite[Theorem~3.7]{Sch:pri} is known as an {\em easy} 
(properly, {\em descent theory}) part of a full structure theorem \cite[Theorem~I]{Sch:pri} for 
Hopf-Galois extensions, recently beautifully extended to a class of coalgebra-Galois extensions 
in \cite{SchSch:Gal}. The difficult part of \cite[Theorem~I]{Sch:pri} involves showing that, when 
appropiate assumptions are made, the bijectivity of the canonical Galois map follows from  its 
surjectivity. In Section~\ref{sec.simple} we show that also in the case of a simple Galois comodule 
of a coring, the surjectivity of the canonical map implies injectivity. In Section~\ref{sec.principal} we  concentrate on Galois comodules which are projective over their 
endomorphism rings. We term such Galois comodules {\em principal comodules}. The 
interest in such comodules stems from non-commutative geometry, in particular from the 
theory of {\em strong connections} \cite{Haj:str} in coalgebra-Galois extensions understood 
as non-commutative principal bundles \cite{BrzMaj:coa}. A certain class of such extensions, 
identified and systematised in \cite{BrzHaj:che} as {\em principal extensions} has 
non-commutative vector bundles, understood as finitely generated projective modules 
(cf.\ \cite{Con:non}), as their associated fibre bundles.  Principal comodules defined 
in the present paper seem to provide a suitable general framework for principal extensions. We derive a full characterisation of principal comodules in the case when the background ring is a field. This is related to the (split-)surjectivity of certian lifting of the canonical map (and hence again resembles the difficult part of Schneider's theorem). We then use this description to prove that, over a field, any entwining structure consisting of an algebra $A$, a coseparable coalgebra $C$ and a bijective entwining map $\psi$ together with a group-like element  in $C$ give rise to a principal extension provided the lifted canonical map is surjective. This can be understood as the entwining structure version of the difficult part of Schneider's Theorem I, and extends recent  theorem of Schauenburg and Schneider  \cite[Theorem~2.5.7]{SchSch:Gal} formulated for a class of Doi-Koppinen entwinings. It also means geometrically that in this case a freeness of the group action induces existence of a (strong) connection on the corresponding principal fibre bundle.

Second class of problems addressed in this paper involves questions, what properties of 
comodules are preserved by morphisms of corings. More precisely, any morphism of corings 
$\cC$ to $\cD$ induces a $\cD$-comodule from a $\cC$-comodule. If a $\cC$-comodule is a 
Galois comodule, is the induced $\cD$-comodule also a Galois $\cD$-comodule? 
Thus in Section~\ref{sec.induce} we determine which morphisms of corings induce principal 
comodules from principal comodules. 
The importance of this induction procedure of principal comodules, and, in particular, principal extensions, in non-commutative geometry 
has been confirmed in recent work \cite{BonCic:bij} in which it has 
been shown that the non-commutative 4-sphere and the corresponding instanton bundle 
constructed in \cite{BonCic:ins} arise from a principal extension of this type. Furthermore, the functor inducing $\cD$-comodules from $\cC$-comodules features prominently in the Kontsevich-Rosenberg approach to non-commutative algebraic geometry \cite{KonRos:smo}, where it is understood as a pull-back of quasi-coherent sheaves over non-commutative stacks, while principal comodules are examples of covers of non-commutative spaces.

The third class of problems discussed in this paper is concerned with duality properties of 
Galois comodules and with the relative injectivity of Galois comodules. In Section~6,  we study 
modules (of the endomorphism ring of a Galois comodule) associated to Galois comodules by 
applying the Hom-functor. The motivation of this construction comes from non-commutative 
geometry, where modules of this kind are understood as fibre bundles associated to 
non-commutative (coalgebra) principal bundles (cf.\ \cite{Brz:mod}). We reveal a remarkable 
duality with respect to the change of arguments in the Hom-bifunctor. This  can again be 
understood in geometric terms as the (generalisation of the) identification of  sections of a 
vector bundle with tensorial zero-forms (functions of type $\rho$). We describe a sufficient 
condition on a Galois comodule $M$ of an $A$-coring $\cC$ that makes any of these associated 
modules a finitely generated projective module (provided a ``fibre'' comodule is a finitely generated projective $A$-module). Finally in Section~7, we connect the relative injectivity of a Galois 
comodule with the properties of the inclusion of the comodule endomorphism rings into the 
module endomorphism ring. This connection then leads to a criterion for faithful flatness of a 
Galois comodule that generalises the criterion introduced for Hopf-Galois extensions  in \cite[Theorem~2.11]{DoiTak:Hop} and for coalgebra-Galois extensions in 
\cite[Proposition~4.4]{Brz:mod}. We also show that if the extension of endomorphism rings of a 
principal comodule is a split extension, then any module,  associated in the way discussed in 
Section~6,   is a finitely generated projective module over the endomorphism ring of the principal comodule (provided a  ``fibre'' comodule is a finitely generated projective  $A$-module).

  \section{Review of corings  and the Galois comodule structure theorem} 
We work over a commutative ring $k$ with a unit. All algebras are over $k$, associative and with a unit. All coalgebras are over $k$, coassociative and with a counit. In a coalgebra $C$, the coproduct is denoted by $\Delta_C$ and the counit by $\eps_C$. The identity morphism for an object $V$ is also denoted by $V$.
 For a ring ($k$-algebra) $R$, the category of right $R$-modules and right $R$-linear maps is denoted by $\M_R$. Symmetric notation is used for left modules. As is customary, we often write $M_R$ to indicate that $M$ is a right $R$-module, etc. When needed, a right action of $R$ on $M_R$ is denoted by $\varrho_M$ and the left action of $R$ on ${}_RN$ is denoted by ${}_N\varrho$. On elements, the actions are denoted by juxtaposition. The dual module of $M_R$ is denoted by $M^*$, while the dual of ${}_RN$ is denoted by ${}^*N$. The product in the endomorphism ring of a right module (comodule) is given by composition of maps, while the product in the endomorphism ring of a left module (comodule) is given by  opposite composition (we always write argument to the right of a function).

Let $A$ be an algebra. A coproduct in an $A$-coring $\cC$ is denoted by  $\DC :\cC\to \cC\ot_A\cC$, and  the counit is denoted by $\eC:\cC\to A$. To
indicate the action of $\DC$ we use the Sweedler sigma notation,
i.e., for all $c\in \cC$,
$$
\DC(c) = \sum c\sw1\ot c\sw 2, \qquad (\DC\ot_A\cC)\circ\DC(c) =
(\cC\ot_A\DC)\circ\DC(c) = \sum c\sw1\ot c\sw 2\ot c\sw 3,
$$
etc. Calligraphic capital letters always denote corings. The category of right $\cC$-comodules and right $\cC$-colinear maps is denoted by $\M^\cC$. Recall that $\M^\cC$ is built upon the category of right $A$-modules, in the sense that there is a forgetful functor $\M^\cC\to \M_A$. In particular, any right $\cC$-comodule is also a right $A$-module, and any right $\cC$-comodule map is right $A$-linear. For a right $\cC$-comodule $M$, $\varrho^M:M\to M\ot_A\cC$ denotes a coaction, and $\Rhom \cC MN$ is the $k$-module of $\cC$-colinear maps $M\to N$. On elements $\varrho^M$ is denoted by the Sweedler notation $\varrho^M(m) = \sum m\sw 0\ot m\sw 1$. Symmetric notation is used for left $\cC$-comodules. In particular, the coaction of a left $\cC$-comodule $N$ is denoted by ${}^N\varrho$, and, on elements, by ${}^N\varrho(n) = \sum n\sw{-1}\ot n\sw 0\in \cC\ot_A N$. Of course, coalgebras are examples of corings, hence the same rules of notation for comodules over a coalgebra as those for comodules over a coring apply. A detailed account of the theory of corings and comodules can be found in \cite{BrzWis:cor}. 

Given a right $\cC$-comodule $M$ and a left $\cC$-comodule $N$ one defines a {\em cotensor product} $M\Box_\cC N$ by the following exact  sequence of $k$-modules:
$$ \xymatrix{0\ar[r]& M\Box_\cC N\ar[r]& M\ot_A N\ar[r]^{\omega_{M,N}\quad}&
             M\ot_A  \cC\ot_A  N,}$$
 where $\omega_{M,N} =\roM \ut_A N - M
\ut_A \Nro$ (i.e., $M\Box_\cC N$ is an equaliser of $\roM \ut_A N$ and  $M
\ut_A \Nro$, where $\roM$ and $\Nro$ are coactions). $M\Box_\cC N$ is a left 
  $S$-module of the endomorphism ring $S=\Rend \cC M$, via $s(\sum_i m_i\ot n_i) = \sum_i s(m_i)\ot n_i$, for all $s\in S$, $\sum_im_i\ot n_i\in M\Box_\cC N$.

Let $\cC$ be an $A$-coring and let $\cD$ be a $B$-coring. A morphism of corings is a pair $(\alpha,\gamma)$, where $\alpha : A\to B$ is an algebra map and $\gamma:\cC\to \cD$ is an $(A,A)$-bimodule map such that 
$$
\chi\circ(\gamma\otimes_A\gamma)\circ\DC = 
\DD\circ\gamma, \quad \eD\circ\gamma = \alpha\circ\eC,
$$
where $\chi: \cD\otimes_A\cD\to \cD\otimes_B\cD$ is the canonical
morphism of $(A,A)$-bimodules induced by $\alpha$. The $(A,A)$-bimodule structure of $\cD$ is induced from the $(B,B)$-bimodule structure via the map $\alpha$ (i.e., $ada' = \alpha(a)d\alpha(a')$, for all $a,a'\in A$ and $d\in \cD$). Such a morphism of corings is explicitly denoted by $(\gamma:\alpha) : (\cC:A)\to (\cD:B)$. In this case any right  $\cC$-comodule $M$ gives rise to a right $\cD$-comodule $M\ot_A B$  with the coaction 
$$
\varrho^{M\ot_A B}: M\ot_A B\to M\ot_A B\ot_B\cD \simeq M\ot_A\cD, \quad  m\ot b \mapsto  \sum m\sw 0\ot\gamma(m\sw 1)b.
$$
For any right $\cC$-comodule $M$ and a right $\cD$-comodule $N$, $\Rhom \cD {M\otimes_AB} N $ is a 
  right $S$-module of the endomorphism ring $S=\Rend \cC M$ via $fs(m\ot b) = f(s(m)\ot b)$, for all $s\in S$, $f\in \Rhom \cD {M\otimes_AB} N $, $m\in M$ and $b\in B$.

 Symmetrically, any left $\cC$-comodule $N$ gives rise to a left $\cD$-comodule $B\ot_A N$. In particular $B\ot_A\cC$ is a left $\cD$-comodule with the coaction 
$${}^{B\ot_A\cC}\varrho :B\ot_A \cC \to \cD\ot_BB\ot_A\cC\simeq \cD\ot_A\cC, \quad b\ot c \mapsto \sum b\gamma(c\sw 1)\ot c\sw 2.
$$ 
Thus for any morphism of corings $(\gamma:\alpha) : (\cC:A)\to (\cD:B)$, there is an associated pair of functors
$$
-\ot_A B : \M^\cC\to \M^\cD, \quad -\Box_\cD(B\ot_A\cC):\M^\cD \to \M_A.
$$
If $\cC$ is flat as a left $A$-module, one shows that $-\Box_\cD(B\ot_A\cC):\M^\cD \to \M^\cC$, and it is a right adjoint of $-\ot_AB$.  We refer to \cite[Section~24]{BrzWis:cor} for more details about morphisms of corings and associated functors. By an {\em $A$-coring morphism}, a morphism of corings $(\gamma:\alpha) : (\cC:A)\to (\cD:A)$ is meant in which $\alpha$ is the identity map (so that only $\gamma$ needs to be specified).

Given right $\cC$-comodules $M$ and  $N$, the $k$-module $\Rhom\cC MN$ is a right module of the endomorphism ring $S=\Rend\cC M$ with the standard action  $fs = f\circ s$, for all $f\in \Rhom \cC MN$, $s\in S$. This defines a functor $\Rhom \cC M- :\M^\cC\to \M_S$. The functor $\Rhom\cC M-$ has the left adjoint $-\ot_SM:\M_S\to \M^\cC$, where, for any $X\in \M_S$, $X\ot_SM$ is a right $\cC$-comodule with the coaction $X\ot_S\roM$. The counit of the adjunction is given by the evaluation map
$$
\varphi_{N}: \Rhom \cC M N \ot_SM\to N, \quad f\ot m\mapsto f(m),
$$
while the unit is $\nu_X: X\to \Rhom \cC M {X\ot_SM}$, $x\mapsto [m\mapsto x\ot m]$ (cf.\ \cite[18.21]{BrzWis:cor}).

This paper is concerned with a special class of comodules introduced in \cite{KaoGom:com} and known as {\em Galois comodules}. The properties of these comodules reflect properties of the above pair of adjoint functors. Let $\cC$ be an $A$-coring,  $M$ be a right $\cC$-comodule and let $S=\Rend \cC M$.  View $\cC$ as a right $\cC$-comodule with the regular coaction $\DC$. 
$M$ is called a {\em Galois (right) 
 comodule} if $M$ is a finitely generated and
projective right $A$-module, and the evaluation
 map  
 $$ \varphi_\cC: \Rhom \cC M\cC\ot_SM\to \cC, \qquad  f\ot m\mapsto f(m),$$
is an isomorphism of right $\cC$-comodules.

An equivalent definition of Galois comodules is obtained by first noting that $M$ is an $(S,A)$-bimodule and $\Rhom \cC M\cC \simeq M^* = \rhom A M A$ as $(A,S)$-bimodules. If $M_A$ is finitely generated projective, then $M^*\ot_S M$ is an $A$-coring with the coproduct $\Delta_{M^*\ot_S M}(\xi\ot m) =  \sum_i \xi\ot e^i \ot \xi^i\ot m$, where $\{e^i\in M,\xi^i\in M^*\}$ is a dual basis of $M_A$, and with the counit $\eps_{M^*\ot_S M}(\xi\ot m) =   \xi(m)$ (cf.\ \cite{KaoGom:com}). The map  $\varphi_\cC$ reduces to the {\em canonical} $A$-coring morphism
$$
\can_M: M^*\ot_S M\to \cC, \qquad \xi\ot m\mapsto \sum \xi(m\sw0)m\sw 1.
$$
$M$ (with $M_A$ finitely generated projective) is a Galois comodule if and only if the canonical map $\can_M$ is an isomorphism of corings. 

The case in which $A$ is a Galois $\cC$-comodule is of fundamental importance. In this case the coaction $\varrho^A :A\to A\ot_A\cC\simeq \cC$ is fully determined by a group-like element $g = \varrho^A(1)\in\cC$. The endomorphism ring $S=\Rend \cC A$ coincides with the subalgebra of $g$-coinvariants in $A$, i.e., $S = \{s\in A\; |\; sg=gs\}$. Obviously, $A$ is a finitely generated projective right $A$-module, $A^* \simeq A$, and $A\ot_S A$ is the {\em Sweedler $A$-coring}, with coproduct $a\ot a'\mapsto a\ot 1\ot 1\ot a'$ and counit $a\ot a'\mapsto aa'$. The canonical map comes out as
$$
\can_A :A\ot_S A\to \cC, \qquad a\ot a'\mapsto aga'.
$$
Thus $A$ is a Galois comodule if and only if $\cC$ is a {\em Galois coring} with respect to $g$, a notion introduced in \cite{Brz:str}.

Main properties of Galois comodules are contained in the Galois comodule structure theorem, which, in part, was first formulated in \cite[Theorem~3.2]{KaoGom:com}. 

\begin{theorem} {\bf (The Galois comodule structure theorem)}
\label{thm.galoisc}
Let $\cC$ be an $A$-coring and $M$ be a right $\cC$-comodule that 
 is finitely generated and projective as a right $A$-module. Set $S=\Rend \cC M$. 

\begin{zlist}
\item The following are equivalent:
 \begin{blist}
 \item $M$ is a Galois comodule that is flat as a left $S$-module.
 \item $\cC$ is a flat left $A$-module and $M$ is a generator in $\M^\cC$.
\item $\cC$ is a flat left $A$-module and, for any $N\in \M^\cC$, the counit of adjunction $\varphi_{N}: \Rhom \cC M N \ot_SM\to N$
is an isomorphism of right $\cC$-comodules.
 \end{blist}
\item The following are equivalent:
  \begin{blist}
 \item $M$ is a Galois comodule  that is faithfully flat as a left $S$-module. 
 \item $\cC$ is a flat left $A$-module  and $M$ is a projective generator in $\cM^\cC$.
 \item  $\cC$ is a flat left $A$-module  and $\Rhom \cC M - : \cM^\cC \to \cM_S$ is an
        equivalence     
        with the inverse $-\ot_SM:\cM_S\to \cM^{\cC}$.
\end{blist}
\end{zlist}
\end{theorem}

For the proof of this theorem we refer to \cite[18.27]{BrzWis:cor} and only point out that the equivalence (b)$\Lra$(c) in (1)  is a consequence of the description of generators in $\M^\cC$ as static comodules in \cite[18.23]{BrzWis:cor}.

\section{Simple Galois comodules}
\label{sec.simple}
The aim of this section is to prove that, for a simple $\cC$-comodule  $M$, to show that $M$ 
is a Galois comodule suffices it to check whether the map $\varphi_M$ is surjective. Recall that 
an object $M$ in an Abelian category is a {\em simple object} provided every monomorphism 
$N\to M$ is either $0$ or an isomorphism. The following characterisation of simple comodules extends 
a theorem of Takeuchi reported in \cite{MasYan:hop}.
\begin{theorem}\label{thm.simple}
Let $\cC$ be an $A$-coring that is flat as a left $A$-module.  Let $M$ 
be a right $\cC$-comodule and  let $S = \Rend \cC M$ be its endomorphism ring.  Then the following are equivalent:
\begin{zlist}
\item $M$ is a simple comodule, i.e., a simple object in $\M^\cC$.
\item $S$ is a division ring and for any right $\cC$-comodule $N$, the evaluation map
$$
\varphi_N :\Rhom \cC MN\otimes_SM\to N, \quad f\ot m \mapsto f(m),
$$
is a monomorphism in $\M^\cC$. 
\end{zlist}
\end{theorem}
\begin{proof}
Since $\cC$  is flat as a left $A$-module, the category $\M^\cC$ is a Grothendieck category (cf.\ \cite[18.14]{BrzWis:cor}), hence, in particular, it is Abelian. 

(1)$\Ra$(2) If $M$ is a simple comodule, then $S$ is a division ring by the Schur lemma. Thus we 
need only to show that, for any $N\in \M^\cC$, the map $\varphi_N$ has a trivial kernel. Note that any element of $\Rhom \cC M N\ot_S M$ (and hence also of the kernel of $\varphi_N$) can be written as a finite sum $\sum_i f_i\otimes m_i$ with the $f_i\in \Rhom \cC M N$ and $m_i \in M$. Since $S$ is a division ring we can always choose the $f_i$ in such a way that they form a free set in the right $S$-module $\Rhom \cC M N$, and we always choose the $f_i$ in this way. Suppose that a simple tensor $f\otimes m$ is in the kernel of $\varphi_N$, i.e., that $f(m) =0$. This means that $m\in \ker f$. On the other hand $M$ is a simple object so that the inclusion monomorphism $0\to \ker f\to M$ is either $0$ or an isomorphism. In the first case the kernel of $f$ is trivial, hence $m=0$, and therefore $f\otimes m =0$. In the other case every $m\in M$ is in the kernel of $M$, hence $f =0$ and $m\otimes f =0$. Thus the kernel of $\varphi_N$ does not contain any non-trivial simple tensors.

Now assume inductively that any non-trivial element consisting of less than $n$ simple tensors cannot be in the kernel of $\varphi_N$, i.e., that $\sum_{i=1}^{n-1} f_i(m_i) = 0$ implies that $m_1 = m_2 =\ldots = m_{n-1} =0$ (as explained, we choose the $f_i$ in such a way that they form a free set). Suppose to the contrary that there exist non-zero $m_n\in M$ and $f_n\in \Rhom \cC M N$ such that
$$
f_1(m_1) + f_2(m_2) +\ldots + f_n(m_n) =0,
$$ 
and $\{f_1, f_2,\ldots , f_n\}$ is a free set in the right $S$-module $\Rhom \cC M N$. This implies that
$$
f_n(M)\cap (\oplus_{k=1}^{n-1} f_k(M)) \neq 0.
$$
Next observe that $f_n(M) \simeq f_nS\ot_S M \simeq M$ via the isomorphism $f_ns\ot m\mapsto sm$, well-defined because $S$ is a division ring. Since $M$ is a simple object, so is $f_n(M)\simeq M$, thus the above intersection property implies that 
$$
f_n(M)\subset \oplus_{k=1}^{n-1} f_k(M). \eqno{(*)}
$$
 Note that $f_k(M) = f_kS\otimes_S M$, and since every $f_kS$ is a finitely generated free $S$-module, there is the following chain of isomorphisms
$$
\Rhom \cC M{ f_k(M)}\simeq f_kS\otimes_S \Rhom \cC MM = f_kS\otimes_S S \simeq f_kS.
$$
Applying $\Rhom \cC M-$ to the inclusion ($*$) we thus obtain
$$
f_n S\subset \Rhom \cC M{\oplus_{k=1}^{n-1} f_k(M)} = \oplus_{k=1}^{n-1} \Rhom \cC M{f_k(M)} \simeq \oplus_{k=1}^{n-1} f_kS.
$$
This, however, contradicts the assumption that the set $\{f_1, f_2,\ldots , f_n\}$ is $S$-free. Hence $
f_1(m_1) + f_2(m_2) +\ldots + f_n(m_n) =0$ implies that $m_1 = m_2 =\ldots = m_{n} =0$ and $\ker\varphi_N =0$ by induction.

(2)$\Ra$(1) Let $f: J\to M$ be a monomorphism in $\M^\cC$ and let $N ={\rm coker} f$. Consider the following commutative diagram with exact rows and columns
$$
\xymatrix{
 0\ar[r] & S\ot_SM\ar[d]_{\pi\ot_SM}\ar[rr]^{\simeq}&& 
  M \ar[r]\ar[d]^p & 0\\
0\ar[r] &  \Rhom \cC M N\otimes_SM\ar[rr]^{\qquad\varphi_N}&& N\ar[d] &\\
&&& 0,& }
$$
where $p$ is the canonical epimorphism and the map $\pi: S\to \Rhom \cC M N$ is given by $\pi(s)(m) = p(s(m))$. It follows that $\varphi_N$ is an isomorphism, and therefore, $\pi\ot_SM$ is an epimorphism. Since $S$ is a division ring, $M$ is a faithfully flat left $S$-module, hence also $\pi$ is an epimorphism.  Furthermore, since $\pi$ is a right $S$-linear map and $S$ is a division ring,  $\ker \pi = 0$ or $\ker \pi = S$. If $\ker \pi =0$, then $\pi$ is an isomorphism, and so is $p$, thus ${\rm coker} f \simeq M$ and therefore $f$ is a zero map. If $\ker \pi =S$, then $\pi$ is the zero map, so $\Rhom \cC M N =0$, i.e., ${\rm coker} f =0$ and $f$ is an isomorphism. Thus $M$ is a simple object.
\end{proof}

Theorem~\ref{thm.simple} leads to the following description of  simple Galois comodules.

\begin{corollary}\label{cor.simple}
Let $\cC$ be an $A$-coring that is flat as a left $A$-module. Then:
\begin{zlist}
\item Every Galois comodule whose endomorphism ring is a division ring is a simple comodule.
\item Let $M$ be a simple right $\cC$-comodule that is finitely generated and projective as a right $A$-module, and let $S=\Rend\cC M$. Then the following are equivalent:
\begin{blist}
\item $M$ is a Galois comodule.
\item For all $N\in \M^\cC$, the evaluation map $\varphi_N : \Rhom\cC MN\ot_S M\to N$ is surjective.
\item The evaluation map $\varphi_\cC : \Rhom\cC M\cC\ot_S M\to \cC$ is surjective.
\item The canonical map $\can_M: M^*\ot_S M\to \cC$, $\xi\ot m\mapsto \sum \xi(m\sw 0)m\sw 1$ is an epimorphism of $A$-corings.
\end{blist}
\end{zlist}
\end{corollary}
\begin{proof}
(1) If $M$ is a Galois comodule with endomorphism ring $S$ that is a division ring, then  $M$ is a flat left $S$-module, hence $\varphi_N$ is an isomorphism for any $N\in\M^\cC$ by Theorem~\ref{thm.galoisc}. Thus $M$ is a simple object in $\M^\cC$ by Theorem~\ref{thm.simple}.

(2) Since $M$ is a simple comodule, $S$ is a division ring, and every $\cC$-comodule is a flat left $S$-module. Hence the implications (a)$\Ra$(b)$\Ra$(c) are obvious. In view of Theorem~\ref{thm.simple}, the evaluation map $\varphi_\cC$ is injective, so condition (c) implies condition (a). Finally the equivalence (c)$\Leftrightarrow$(d) follows from the hom-tensor isomorphism $\Rhom\cC M \cC \simeq \Rhom \cC M {A\ot_A\cC} \simeq \rhom A M A$.
\end{proof}

Note in passing that, by extracting the key features of  Theorem~\ref{thm.simple} and using a metatheorem of Abelian categories (cf.\ \cite[Chapter~4]{Fre:abe}), one can obtain a characterisation of simple objects in general Abelian categories. The key features are that  the 
functor $-\ot_S M$ is the left adjoint of the functor 
$\Rhom \cC M -$ and that the map $\varphi_N$ is a counit of 
this adjunction. Furthermore, the statement of the theorem is a compound diagrammatic statement. Finally, the facts that $A$ is a (trivial) $A$-coring, and the category of right $A$-comodules is isomorphic to the category of  right $A$-modules assure that Theorem~\ref{thm.simple} holds for any category of modules. Thus a metatheorem of Abelian categories together with the Mitchell Embedding Theorem (cf.\  \cite[Chapter~7]{Fre:abe}) lead to the following characterisation of simple objects. In 
any Abelian category $\C$, an object $M$ such that ${\rm Mor}_\C(M,-)$ has the left 
adjoint is simple if and only if its endomorphism ring is a division ring and, for any 
object, the counit of the adjunction is a monomorphism (compare
 characterisation of adjoints in \cite[Chapitre V]{Gab:cat}). This is probably well-known to category 
theorists (although we were not able to find a reference).

  \section{Principal comodules} 
\label{sec.principal}
In this section we introduce and study 
the following class of Galois comodules.
\begin{definition} A Galois right $\cC$-comodule $M$ is said to be a {\em principal comodule} provided it is a projective left module of its endomorphism ring $S=\Rend \cC M$.
\label{def.principal}
\end{definition}

The prime interest in studying principal comodules stems from
non-commutative geometry. One can argue that a principal comodule is as
close an object as abstractly possible to the notion of a {\em principal
extension} introduced recently in \cite{BrzHaj:che}. The latter is an example of a principal comodule of a coring associated to an {\em entwining structure}.

Recall from \cite{BrzMaj:coa} that an entwining structure $(A,C)_\psi$ consists of of a $k$-algebra $A$, a $k$-coalgebra
$C$ and an $k$-module map 
$ \psi: C\tens_{k} A\to A\tens_{k} C$ rendering commutative 
the
following  {\em bow-tie diagram}                                             
$$
\xymatrix{
& C\ot_k  A\ot_k  A \ar[ddl]_{\psi\otimes_k A} 
\ar[dr]^{C\otimes_k \mu}  &  &             
C\ot_k  C\ot_k  A \ar[ddr]^{C\ot_k\psi}& \\
& & C\ot_k  A \ar[ur]^{\Delta_C\ot_k A} 
\ar[dr]^{\eps_C\ot_k A}\ar[dd]^{\psi}  & &\\
A\ot_k C\ot_k  A \ar[ddr]_{A\ot_k\psi} & C \ar[ur]^{C\otimes_k 
\iota} \ar[dr]_{\iota\otimes_k C}& &A  & C\ot_k  A\ot_k  C 
\ar[ddl]^{\psi\otimes_k C}\\
& &  A\ot_k  C \ar[ur]_{A\otimes_k\eps_C} \ar[dr]_{A\otimes_k\Delta_C} & &\\
& A\ot_k  A\ot_k  C \ar[ur]_{\mu\otimes_k C} &  & A\ot_k C\ot_k  C & , }
$$
where $\mu$ is the product in $A$ and $\iota:k\to A$ is the unit map. The map $\psi$ is known as an {\em entwining map}\index{entwining map}, and $C$ and $A$ are said to be {\em entwined} by $\psi$. As explained in \cite{Brz:str}, given an entwining structure $(A,C)_\psi$, $\cC= A\ot_kC$ is an $A$-coring with $A$-multiplications $a(a'\ot c)a'' = aa'\psi(c\ot a'')$,  coproduct $\DC: A\ot_k C\to A\ot_k C\ot_AA\ot_k C\simeq A\ot_k C\ot_k C$,  $a\ot c \mapsto  a\ot \Delta_C(c)$, and  counit $\eC(a\ot c) = a\eps_C(c)$. For more information about entwining structures and their connection with Hopf-type modules we refer to \cite{CaeMil:gen}.

\begin{example}
\label{ex.principal}
Let $k$ be a field, $C$ be a coalgebra and $A$ an algebra and a right $C$-comodule via
$\varrho^A:A\to A\ot_k C$. Let
$S=A^{co C}:=\{s\in A~|~\varrho^A(sa)=s\varrho^A(a),\ \forall a\in A\}$, denote the subalgebra of $C$-coinvariants of $A$. 
The inclusion of algebras $S\subseteq A$ is called a
 {\em principal $C$-extension}
 iff
\begin{zlist}
\item
$
\can: A\ot_SA{\to} A\ot_k C,\;a\ot a'\mapsto a\varrho^A(a')
$
is bijective (the Galois  condition);
\item 
$A$ is $C$-equivariantly projective as a left $S$-module, i.e., there exists a left $S$-module, right $C$-comodule section of the product $S\ot_k A\to A$ (existence of a strong
connection);
\item  $\psi:C\ot_k A{\to} A\ot_k C$, $c\ot a\mapsto
\can(\can^{-1}(1\ot c)a)$ is bijective (invertibility of the canonical entwining);
\item 
there is a group-like element $e\in C$ such that $\varrho^A(a)=\psi(e\ot a)$,
$\forall a\in
A$ (co-augmentation).
\end{zlist}

If $S\subseteq A$ is a principal $C$-extension, then $A$ is a principal comodule of the coring $\cC = A\otimes_k C \simeq A\ot_S A$.

\end{example}
\begin{proof}
By  \cite[Theorem~2.7]{BrzHaj:coa}, the map $\psi$ entwines $A$ with $C$, so $\cC= A\ot_kC$ is an $A$-coring of a type described above. Furthermore,  $A$ is a right $\cC$-comodule with the coaction $\varrho^A: A\to A\ot_k C \simeq A\ot_A A\ot_k C = A\ot_A\cC$. This means that $g=\varrho^A(1) = 1\ot e$ is a group-like element in $\cC$, and $\varrho^A(a) = ga$ for all $a\in A$.
By the obvious identification $\rend AA  \simeq A$ we obtain
$$
\Rend \cC A = \{ f\in \rend AA \; |\; \forall a\in A, \; gf(a) = f(a)g\} = \{s\in A \: |\; sg=gs\} = S.
$$
The Galois condition implies that $\cC$ is a Galois coring, hence  $\cC \simeq A\ot_S A \simeq A^*\ot_S A$, i.e., $A$ is a Galois comodule. Finally, using the above identification of $\Rend \cC A$, the evaluation map $\Rend \cC A\ot_k A\to A$ becomes simply the product map $S\ot_k A\to A$, $s\ot a\mapsto sa$, and hence the existence of a strong connection implies that $A$ is a principal comodule.
\end{proof}

A principal $C$-extension can be seen as a non-commutative
version of a principal bundle. Since  $k$ is a field (this is the main case of interest from the point of view of non-commutative geometry), the coring $A\ot_kC$ in Example~\ref{ex.principal} is free as a left $A$-module. Then  the bijectivity of the canonical entwining map (condition (3) in Example~\ref{ex.principal}) implies that $A\ot_kC$ is also free as a right $A$-module. Thus one could get even closer to a principal
extension by considering principal comodules over corings that are free
as left and right $A$-modules. For the purpose of a general exposition in this paper, however, this would be
an unnecessary restriction on a coring.

Since a principal comodule $M$  is a projective $S$-module,  it is a flat $S$-module (thus, in particular ${}_A\cC$ is a flat module and $M$ is a generator in $\M^\cC$ by the Galois comodule structure theorem). In fact, the principality of a comodule implies faithful flatness. More precisely one proves the following
\begin{theorem}
Let $M$ be a principal $\cC$-comodule that is faithfully flat as a $k$-module and set $S=\Rend \cC M$. Then $M$ is a faithfully flat left $S$-module.
\label{thm.fflat}
\end{theorem}
\begin{proof} 
Consider an epimorphism $f:V\to W$ of right $\cC$-comodules. Since ${}_SM$ is flat, both evaluation maps $\varphi_V$ and $\varphi_W$ in Theorem~\ref{thm.galoisc}(1)(c) are isomorphisms in $\M^\cC$. Note that
$$
\varphi_W^{-1}\circ f \circ \varphi_V = \Rhom \cC M f\ot_SM,
$$
and thus we obtain an exact sequence of right $\cC$-comodule maps
$$\xymatrix{ 
     \Rhom \cC M V\ot_SM\ar[rrr]^{\Rhom \cC M f\ot_SM}&&& 
   \Rhom \cC M W\ot_SM\ar[r] & 0.}
$$
Let $\sigma: M\to S\ot_k M$ be a left $S$-linear splitting of the $S$-action  ${}_M\varrho$. Then the above exact sequence leads to the following commutative diagram with exact top row and split-exact columns
$$
\xymatrix{0 \ar[d] &&& 0\ar[d] &\\
 \Rhom \cC M V\ot_SM\ar@<1ex>[dd]^{\Rhom \cC M V\ot_S\sigma}\ar[rrr]^{\Rhom \cC M f\ot_SM}&&& 
   \Rhom \cC M W\ot_SM\ar@<1ex>[dd]^{\Rhom \cC M W\ot_S\sigma} \ar[r] & 0\\
&&&&\\
 \Rhom \cC M V\ot_kM\ar@<1ex>[uu]^{\Rhom \cC M V\ot_S{}_M\varrho}\ar[rrr]^{\Rhom \cC M f\ot_kM}&&& 
  \Rhom \cC M W\ot_kM\ar@<1ex>[uu]^{\Rhom \cC M W\ot_S{}_M\varrho} & }
$$
Therefore the map $\Rhom \cC M f\ot_kM$ is surjective (for any $x\in  \Rhom \cC M W\ot_kM$ is of the form $(\Rhom \cC M f\ot_kM)(y)$, where $y$ is such that $(\Rhom \cC M f\ot_SM)(y) = (\Rhom \cC M W\ot_S\,{}_M\!\varrho)(x)$). Since $M$ is a faithfully flat $k$-module, also the map $\Rhom \cC M f : \Rhom \cC M V\to \Rhom \cC M W$ is surjective. This implies that $M$ is  a projective object in $\M^\cC$ (cf.\ \cite[18.20]{BrzWis:cor}). Thus $M$ is a projective generator in $\M^\cC$ and the Galois comodule structure theorem implies that ${}_SM$ is a faithfully flat module. 
\end{proof}

Thus, by the Galois comodule structure theorem, every  principal $\cC$-comodule $M$ that is faithfully flat as a $k$-module is a finitely generated projective generator in $\M^\cC$, and it induces the category equivalence $-\ot_S M: \M_S\to \M^\cC$. 

In  case $k$ is a field, one can derive a description of principal comodules which resembles the difficult part of the Schneider theorem. To facilitate this description, note that in parallel to the theory of Galois right comodules one develops the theory of Galois left comodules. Consider a left comodule $N$ of an $A$-coring $\cC$ with the endomorphism ring $T = \Lend \cC N$. The product in $T$ is given by opposite composition, i.e., $tt' = t'\circ t$. This makes $N$ into a right $T$-module, with multiplication $nt = t(n)$. $N$ is called a {\em Galois (left) comodule} provided it is  finitely generated projective left $A$-module and the evaluation
 map 
 $$ \hphi_\cC: N\ot_T\Lhom\cC N\cC\to \cC, \; n\ot f\mapsto f(n),$$
is an isomorphism of left $\cC$-comodules.  Equivalently, $N$ is a Galois left $\cC$-comodule provided ${}_AN$ is finitely generated projective and the left canonical map
$$
{}_N\can : N\ot_T {}^*N \to \cC,\qquad n\ot \xi\mapsto \sum n\sw{-1}\xi(n\sw 0),
$$
is an isomorphism of $A$-corings.

Any Galois right $\cC$-comodule $M$ gives rise to a Galois left $\cC$-comodule. First, since $M$ is a finitely generated projective right $A$-module, the dual module $M^* = \rhom AMA$ is a left $\cC$-comodule with the coaction determined by
\begin{equation}\label{left.coa1}
\sum \xi(m\sw 0)m\sw 1 = \sum \xi\sw{-1}\xi\sw 0(m), \qquad \forall \xi\in M^*, m\in M.
\end{equation}
Explicitly, in terms of a dual basis $\{e^i\in M,\xi^i\in M^*\}_{i=1,\ldots, n}$ of $M$, the left coaction comes out as
$$
{}^{M^*\!\!}\varrho : M^*\to \cC\ot_A M^*,\qquad \xi\mapsto \sum_i\xi(e^i\sw 0)e^i\sw 1\ot \xi^i.
$$
This definition of the left $\cC$-coaction also implies the following equality for the elements of a dual basis
\begin{equation}\label{left.coa.2}
\sum_i e^i\sw 0\ot e^i\sw 1\ot \xi^i = \sum_i e^i\ot \xi^i\sw{-1}\ot \xi^i\sw 0.
\end{equation}
Second, there is a ring isomorphism $\Gamma_M: S=\Rend \cC M\to \Lend \cC {M^*}$. Explicitly, for all $s\in S$, $\xi\in M^*$, $t\in \Lend \cC {M^*}$ and $m\in M$, the map $\Gamma_M$ and its inverse $\Gamma_M^{-1}$ are  given by 
\begin{equation}\label{Gamma}
\Gamma_M(s)(\xi) = \sum_i\xi(s(e^i))\xi^i,\qquad \Gamma_M^{-1}(t)(m) = \sum_i e^it(\xi^i)(m).
\end{equation}
Thus $M^*$ is a right $S$-module with the multiplication $\xi s = \Gamma_M(s)(\xi) = \xi\circ s$. Note that $\Gamma_M$ is a restriction of the canonical isomorphism $\rend A M\to \lend A{M^*}$, which we also denote by $\Gamma_M$. Third, with this isomorphism of endomorphism rings, and with the help of equation (\ref{left.coa1}), the left canonical map comes out as ${}_{M^*}\can : M^*\ot_S M\to \cC$, ${}_{M^*}\can (\xi\ot m) = \sum \xi\sw{-1} \xi\sw 0(m) = \sum\xi(m\sw 0)m\sw 1$ and thus ${}_{M^*}\can = \can_M$. Since $\can_M$ is an isomorphism of corings, so is ${}_{M^*}\can$, and therefore $M^*$ is a Galois left $\cC$-comodule. 

This interlude on Galois left comodules allows one to state the following theorem which is a coring origin of \cite[Theorem~2.3.4]{SchSch:Gal}.
\begin{theorem}\label{thm.princ.field}
Let $k$ be a field, $\cC$ an $A$-coring and $M$ a right $\cC$-comodule that is finitely generated and projective as a right $A$-module. Set $S=\Rend \cC M$. View  $M^*\ot_kM$ as a left $\cC$-comodule via ${}^{M^*\!\!}\varrho\ot_k M$. Then the following statements are equivalent:
\begin{zlist}
\item $\cC$ is a flat left $A$-module and the map
$$
\tcan_M :M^*\ot_k M\to \cC, \qquad \xi\ot m\mapsto\sum \xi(m\sw 0)m\sw 1,
$$
is a split epimorphism of left $\cC$-comodules.
\item $M$ is a principal right $\cC$-comodule.
\end{zlist}
\end{theorem}

The proof of Theorem~\ref{thm.princ.field} makes use of the following
\begin{lemma}\label{lemma.princ.field}
Let $k$ be a field, $\cC$ an $A$-coring and $N$ a left $\cC$-comodule that is finitely generated and projective as a left $A$-module. Then for any $k$-vector space $V$,
$$
\Lhom \cC N {N\ot_k V}\simeq \Lend \cC N\ot_k V,
$$
where $N\ot_k V$ is a left $\cC$-comodule with the coaction ${}^{N\!\!}\varrho\ot_k V$.
\end{lemma}
\begin{proof}
Since ${}_AN$ is finitely generated and projective, there is a $k$-linear isomorphism $\theta:  \lend A N\ot_k V\to \lhom A N {N\ot_k V}$. Explicitly, for all $t\in \lend AN$, $v\in V$ and $n\in N$, $\theta(t\ot v)(n) = t(n)\ot v$. Clearly, if $t\in \Lend\cC N$, then $\theta(t\ot v)\in \Lhom \cC N {N\ot_k V}$. Thus we only need to check whether the inverse $\theta^{-1}$ of $\theta$ restricted to $\Lhom \cC N {N\ot_k V}$ has the image in $\Lend \cC N\ot_k V$. 

To write $\theta^{-1}$ explicitly, choose a dual basis $\{\xi^i\in N, \;e^i\in {}^*N\}_{i=1,\ldots l}$ of ${}_AN$ and, for all $f\in \Lhom \cC N {N\ot_k V}$ and $n\in N$, write $f(n) = \sum f(n)\suc1\ot f(n)\suc2\in N\ot_kV$. Then
$$
\theta^{-1}(f) = \sum_i e^i(-)f(\xi^i)\suc1\ot f(\xi^i)\suc2.
$$
Since $f$ is a left $\cC$-comodule map, for all $n\in N$,
\begin{eqnarray*}
\sum n\sw{-1}\ot\theta^{-1}(f)(n\sw 0) &=& \sum_i n\sw{-1}\ot e^i(n\sw 0)f(\xi^i)\suc1\ot f(\xi^i)\suc2\\
&=& \sum_i n\sw{-1}\ot f(e^i(n\sw 0)\xi^i)\suc1\ot f(e^i(n\sw 0)\xi^i)\suc2\\
&=&  \sum n\sw{-1}\ot f(n\sw 0)\suc1\ot f(n\sw 0)\suc2\\
&=& \sum f(n)\suc1 \sw{-1}\ot f(n)\suc1\sw 0\ot f(n)\suc2.
\end{eqnarray*}
The second equality follows from the left $A$-linearity of $f$, the third one is a consequence of  the properties of a dual basis, and the last equality results from $\cC$-colinearity of $f$. On the other hand, again using the properties of a dual basis and the $A$-linearity of $f$, we obtain
\begin{eqnarray*}
({}^{N\!\!}\varrho\ot_k V)(\theta^{-1}(f)(n)) &=& \sum_i (e^i(n)f(\xi^i)\suc1)\sw{-1}\ot (e^i(n)f(\xi^i)\suc 1)\sw 0 \ot f(\xi^i)\suc 2\\
&=& \sum_i f(e^i(n)\xi^i)\suc 1\sw{-1}\ot f(e^i(n)\xi^i)\suc 1\sw 0 \ot f(e^i(n)\xi^i)\suc 2\\
&=& \sum f(n)\suc 1 \sw{-1}\ot f(n)\suc 1\sw 0\ot f(n)\suc 2.
\end{eqnarray*}
Thus $\im (\theta^{-1}\mid_{\Lhom \cC N {N\ot_k V}})\subseteq \Lend \cC N\ot_k V$, hence $\theta$ restricts to an isomorphism $\Lend \cC N\ot_k V \simeq\Lhom \cC N {N\ot_k V}$, as required.
\end{proof}

{\sl Proof of Theorem~\ref{thm.princ.field}.} (1) $\Ra$ (2) Lemma~\ref{lemma.princ.field} and the identification of left $\cC$-endomorphisms of $M^*$ with $S$ leads to the isomorphism $\Lhom \cC {M^*}{M^*\ot_k M}\simeq \Lend \cC {M^*} \ot_k M \simeq S\ot_k M$. Furthermore, $\Lhom \cC {M^*}\cC\simeq {}^*(M^*)\simeq M$. Since $\cC$ is a direct summand of $M^*\ot_k M$ as a left $\cC$-comodule, and $\Lhom \cC {M^*}-$ is a functor from ${}^\cC\M$ to ${}_S\M$, the above isomorphisms imply that $M$ is a direct summand of $S\ot_kM$ as a left $S$-module. Since $k$ is assumed to be a field, $S\ot_kM$ is a free $S$-module, hence $M$ is a projective left $S$-module.

The first step in the proof that $M$ is a Galois comodule is to consider the following commutative diagram
$$
\xymatrix{
M^*\ot_S \Lhom \cC {M^*} {M^*\ot_k M}\ar[dr]^{\simeq}\ar[rr]^{\hphi_{M^*\ot_k M}}&& 
  M^*\ot_k M\\
& M^*\ot_SS\ot_k M \ar[ur]^{\simeq} &, }
$$
where the first isomorphism follows from Lemma~\ref{lemma.princ.field} and the discussion above. Thus $\hphi_{M^*\ot_k M}$ is an isomorphism. Next we can consider the following diagram, which is commutative in all possible directions since $\hphi$ is a natural transformation,
$$
\xymatrix{M^*\ot_S \Lhom \cC {M^*} {M^*\ot_k M}\ar@<1ex>[dd]^{M^*\ot_S \Lhom \cC {M^*} {\tcan_M}}\ar[rrr]^{\hphi_{M^*\ot_k M}}&&& 
   M^*\ot_k M\ar@<1ex>[dd]^{\tcan_M}\\
&&&\\
M^*\ot_S \Lhom \cC {M^*} \cC\ar@<1ex>[uu]\ar[rrr]^{\hphi_\cC}&&& 
\cC\ar@<1ex>[uu] .}
$$
The upward pointing arrows are sections of $M^*\ot_S \Lhom \cC {M^*} {\tcan_M}$ and $\tcan_M$ respectively. Since $\hphi_{M^*\ot_k M}$ is an isomorphism, the map $\hphi_\cC$ is one-to-one and onto (it is a $k$-linear isomorphism). With the help of identifications $\Lhom \cC {M^*}\cC\simeq M$ and $M^*\simeq \Rhom \cC M\cC$ we can construct yet another commutative diagram
$$
\xymatrix{
M^*\ot_S \Lhom \cC {M^*}\cC\ar[dr]^{\simeq}\ar[rr]^{\hphi_\cC}&& 
 \cC\\
& \Rhom \cC M\cC \ot_SM \ar[ur]^{\varphi_\cC} &. }
$$
Since $\hphi_\cC$ is one-to-one and onto, so is $\varphi_\cC$. By assumption, $\cC$ is a flat left $A$-module, so $\M^\cC$ is an Abelian category. Since $\varphi_\cC$ is a one-to-one and onto morphism in $\M^\cC$, it is an isomorphism. Thus $M$ is a Galois right $\cC$-comodule that is projective over its endomorphism ring, i.e., $M$ is a principal comodule.

(2) $\Ra$ (1) Applying functor $M^*\ot_S - : {}_S\M\to {}^\cC\M$ to a left $S$-linear section of the multiplication $S\ot_k M\to M$, one obtains a left $\cC$-comodule section of the canonical surjection $M^*\ot_kM\to M^*\ot_S M$. Composed with an $A$-coring (hence also a left $\cC$-comodule) map $\can^{-1}_M:\cC\to M^*\ot_SM$ this gives the required section of $\tcan_M$.
 \endproof
 
 Theorem~\ref{thm.princ.field} leads to the main geometric result of this section, namely to a condition for an entwining structure to give rise to a principal extension or a non-commutative principal fibre bundle. First recall that a coalgebra $C$ is said to be {\em coseparable} provided  
the coproduct has a retraction in the category of $C$-bicomodules.
Equivalently, $C$ is a coseparable coalgebra if there exists a {\em
cointegral}, i.e., a $k$-module map $\delta :C\ot_kC\to k$ such
that $\delta\circ\Delta_C = \eps_C$ and, for all $c,c'\in C$,
$$
\sum c\sw 1\delta(c\sw 2\ot c') = \sum \delta(c\ot c'\sw 1)c'\sw 2.
$$
This equality is known as the  {\em colinearity} of the cointegral $\delta$. Any cosemisimple coalgebra over an algebraically closed field is coseparable (cf.\ \cite[Remark~2.5.3]{SchSch:Gal}).

\begin{theorem}\label{thm.main}
Let $k$ be a field and $(A,C)_\psi$ an entwining structure such that the map $\psi$ is bijective. Suppose that $e\in C$ is a group-like element and view $A$ as a right $C$-comodule with the coaction $\varrho^A:A\to A\ot_kC$, $a\mapsto \psi(e\ot a)$. If $C$ is a coseparable coalgebra and the (lifted) canonical map
$$
\tcan: A\ot_kA{\to} A\ot_k C,\;a\ot a'\mapsto a\varrho^A(a')
$$
is surjective, then $A$ is a principal $C$-extension of the coinvariants 
$S=A^{co C}$.
\end{theorem}
\begin{proof}
The strategy for the proof is first to use  Theorem~\ref{thm.princ.field} (with $M=A$) to show that $A$ is a principal comodule of the coring $\cC = A\ot_kC$ corresponding to $(A,C)_\psi$, and then to show that the projectivity of ${}_SA$ implies the $C$-equivariant projectivity (cf.\ \cite[Proposition~2.5.4]{SchSch:Gal}).

Following Theorem~\ref{thm.princ.field} we need to construct a left $\cC$-comodule splitting of $\tcan$. The right $C$-coaction $\varrho^A$ can be understood as a right $\cC$-coaction corresponding to a group-like element $g=1\ot e\in \cC$, i.e., $\varrho^A(a) = \psi(e\ot a) = (1\ot e) a = ga$ (cf.\ proof of Example~\ref{ex.principal}). Thus the left $\cC$-coaction of $A^*\simeq A$ comes out as $a\mapsto ag = a\ot e$. Since $\psi$ is bijective, this left $\cC$-coaction  gives rise to a left $C$-coaction $\Aro(a) = \psi^{-1}(a\ot e) :=\sum a\sw{-1}\ot a\sw 0$.  In view of the isomorphism $\lhom A {A\ot_k C}{A\ot_k A} \simeq {\rm Hom}_k(C, A\ot_k A)$, any left $\cC$-comodule map $f:A\ot_kC\to A\ot_k A$ can be identified with a $k$-linear map $\hat{f}: C\to A\ot_kA$ such that, for all $c\in C$, writing $\hat{f}(c) := \sum \hat{f}(c)\suc 1\ot \hat{f}(c)\suc 2\in A\ot_k A$,
$$
\sum \psi(c\sw 1\ot  \hat{f}(c\sw 2)\suc 1)\ot  \hat{f}(c)\suc 2 = \sum  \hat{f}(c)\suc 1\ot e\ot  \hat{f}(c)\suc  2\in A\ot_k C\ot_k A.
$$
Applying $\psi^{-1}\ot_k A$ we thus, equivalently, obtain the condition
$$
(C\ot_k \hat{f})\circ \Delta_C = (\Aro \ot_k A)\circ\hat{f}. \eqno{(*)}
$$
Since $\tcan$ is surjective, it has a $k$-linear section $\tau: A\ot_k C\to A\ot_k A$. Define $\hat{\tau} : C\to  A\ot_k A$ by $\hat{\tau}(c) = \tau(1\ot c)$.  Let $\delta$ be a cointegral of $C$ and define a $k$-linear map $\hat{\kappa}: C\to A\ot_k A$, by
$$
\hat{\kappa} = (\delta\ot_k A\ot_k A)\circ (C\ot_k\Aro\ot_k A)\circ(C\ot_k\hat{\tau})\circ \Delta_C.
$$
Using the colinearity of $\delta$ one easily checks that $\hat{\kappa}$ has the property ($*$). Therefore the map $\kappa: A\ot_k C\to A\ot_k A$, 
$\kappa(a\ot c) = a\hat{\kappa}(c)$
is a left $\cC$-comodule morphism. We aim to prove that $\kappa$ is a section of $\tcan$.

To this end, first introduce the {\em $\alpha$-notation} for an entwining map and its inverse, i.e., for all $a\in A$, $c\in C$, write
$$\psi(c\ot a) = \sum_\alpha a_\alpha\ot c^\alpha , \qquad \psi^{-1}(a\ot c) = \sum_A c_A\ot a^A.
$$
In this notation the left and right pentagon conditions in the bow-tie diagram read respectively,
$$
\psi(c\ot aa') = 
\sum_{\alpha,\beta} 
     a_\alpha a'_\beta \ot
 c^{\alpha\beta}, \qquad \sum_\alpha  a_\alpha\ot c^\alpha\sw 1\ut c^\alpha\sw 2 =
  \sum_{\alpha,\beta}a_{\beta\alpha}\ot 
c\sw 1^\alpha\ut c\sw 2^\beta .
  $$
Since $\psi^{-1}$ is the inverse of $\psi$,
$$
a\ot c = \sum_{\alpha,A} a_\alpha^A \ot c_A^\alpha.
$$
Second, write $\hat{\tau}(c) = \sum c\suc 1\ot c\suc 2\in A\ot_kA$, so that the map $\kappa$ explicitly reads
$$
{\kappa} (a\ot c) = \sum a\delta(c\sw 1\ot c\sw 2\suc 1\sw{-1})c\sw 2\suc 1\sw 0 \ot c\sw 2\suc 2.
$$
Note that, since $\hat{\tau}$ is obtained from a $k$-linear section of $\tcan$, for all $c\in C$,
$$
\sum c\suc 1c\suc 2\sw 0\ot c\suc 2\sw 1  = 1\ot c. \eqno{(**)}
$$
In view of the definitions of left and right $C$-coactions, for all $c\in C$,
$$
\sum c\suc 1\sw{-1}\ot c\suc 1\sw 0c\suc 2\sw 0 \ot c\suc 2\sw 1 = \sum \psi^{-1}(c\suc 1\ot e) \psi(e\ot c\suc 2) = \sum_{\alpha, A} e_A\ot c\suc 1^Ac\suc 2_\alpha \ot e^\alpha.
$$
Apply $\psi\ot C$ to this identity and compute
\begin{eqnarray*}
\sum \psi( c\suc 1\sw{-1}\ot c\suc 1\sw 0c\suc 2\sw 0) \ot c\suc 2\sw 1&=& \sum_{\alpha, A}\psi( e_A\ot c\suc 1^Ac\suc 2_\alpha) \ot e^\alpha\\
&=&  \sum_{\alpha,\beta,\gamma, A}  c\suc 1^A_\gamma c\suc 2_{\alpha\beta}\ot e_A^{\gamma\beta} \ot e^\alpha\\
&=&  \sum_{\alpha,\beta}  c\suc 1 c\suc 2_{\alpha\beta}\ot e^{\beta} \ot e^\alpha\\
&=&  \sum_{\alpha}  c\suc 1 c\suc 2_{\alpha}\ot e^{\alpha}\sw 1 \ot e^\alpha\sw 2\\
&=&  \sum  c\suc 1 c\suc 2\sw 0\ot c\suc 2\sw 1 \ot c\suc 2\sw 2\\
&=& \sum 1\ot c\sw 1\ot c\sw 2.
\end{eqnarray*}
The second equality follows from the left pentagon in the bow-tie diagram, while the third one is the consequence of the fact that $\psi^{-1}$ is the inverse of $\psi$. Next, observing that $e$ is a group-like element and using the right pentagon in the bow-tie diagram we arrive at the fourth equality. The remaining two equalities follow from the definition of the right $C$-coaction and from equation ($**$). The left triangle in the bow-tie diagram implies that $\psi^{-1}(1\ot c) = c\ot 1$, thus applying $\psi^{-1}$ to the equality just derived, we conclude that
$$
\sum c\suc 1\sw{-1}\ot c\suc 1\sw 0c\suc 2\sw 0 \ot c\suc 2\sw 1 = \sum c\sw 1\ot 1\ot c\sw 2. \eqno{(***)}
$$
Now equation ($***$) facilitates the following computation
\begin{eqnarray*}
(\tcan\circ\kappa)(a\ot c) &=& \sum a\delta(c\sw 1\ot c\sw 2\suc 1\sw{-1}) c\sw 2\suc 1\sw 0c\sw 2\suc 2\sw 0 \ot c\sw 2\suc 2\sw 1\\
&=& \sum a\delta(c\sw 1\ot c\sw 2)\ot c\sw 3 = \sum a\ot \eps_C(c\sw 1)c\sw 2 = a\ot c,
\end{eqnarray*}
where the penultimate equality is the consequence of the fact that $\delta$ is a cointegral. Thus we have proven that $\kappa$ is a left $\cC$-colinear section of the canonical map $\tcan$ so that $A$ is a principal right $\cC$-comodule by Theorem~\ref{thm.princ.field}. This means in particular that $A$ is a Galois $C$-extension of $S$, i.e., condition (1) in Example~\ref{ex.principal} is satisfied. Furthermore, the uniqueness of the canonical entwining map (cf.\  \cite[Theorem~2.7]{BrzHaj:coa}) implies that $\psi$ is the canonical entwining map, hence it is bijective as required for condition (3) in Example~\ref{ex.principal}. Obviously, condition   (4) in Example~\ref{ex.principal} is also satisfied. Thus we only need to construct a section satisfying condition (2) in Example~\ref{ex.principal}.

Since $A$ is a principal $\cC$-comodule it is a projective left $S$-module, hence there is a left $S$-module section $\tilde{\sigma}: A\to S\ot_k A$ of the product. Using the cointegral $\delta$, construct the map
$$
\sigma: A\to S\ot_k A, \quad \sigma = (S\ot_k A\ot_k \delta)\circ(S\ot_k \varrho^A\ot_k C) \circ(\tilde{\sigma}\ot_k C)\circ\varrho^A.
$$
Since $\sigma$ is a composition of $S$-module maps it is a left $S$-module map. Using the colinearity  of a cointegral,  one easily shows that $\sigma$ is a right $C$-colinear  map. To show that $\sigma$ is a section of the product map,  we denote $\tilde{\sigma}(a) = \sum a\suc 1\ot a\suc 2\in S\ot_k A$, so that the map $\sigma$ explicitly reads,
$$
\sigma(a) = \sum a\sw 0\suc 1\ot a\sw 0\suc 2\sw 0\delta(a\sw 0\suc 2\sw 1\ot a\sw1).
$$
Remember that $\tilde{\sigma}$ is a section of the product map ${}_A\varrho: S\ot_k A\to A$, $s\ot a\to sa$, so that for all $a\in A$, $\sum a\suc 1a\suc 2 = a$. Remember also that ${\varrho}^A$ is a left $S$-linear map ($S$ is  the $\cC$-endomorphism ring of $A$!). With these facts at hand we can compute
\begin{eqnarray*}
({}_A\varrho\circ \sigma)(a) &=& \sum a\sw 0\suc 1 a\sw 0\suc 2\sw 0\delta(a\sw 0\suc 2\sw 1\ot a\sw1)\\
&=& \sum (a\sw 0\suc 1 a\sw 0\suc 2)\sw 0\delta((a\sw 0\suc 1 a\sw 0\suc 2)\sw 1\ot a\sw1)\\
&=& \sum a\sw 0\delta(a\sw 1\ot a\sw 2) =a.
\end{eqnarray*}
Therefore, $\sigma$ is a right $C$-colinear, left $S$-linear section of the product map, so that also condition (2) in Example~\ref{ex.principal} is satisfied, hence $A$ is a principal $C$-extension of $S$.
\end{proof}

Since any comodule of a coseparable coalgebra over a field is an injective comodule, Theorem~\ref{thm.main} can be understood as an entwining structure version of the `difficult part' of Schneider's structure theorem \cite[Theorem~I]{Sch:pri}. As a special case one obtains 
\begin{corollary}{\cite[Theorem~2.5.7]{SchSch:Gal}}
Let $k$ be a field, $H$ be a Hopf algebra with a bijective antipode and let $A$ be a right $H$-comodule algebra with a coaction $\bar{\varrho}^A:A\to A\ot_k H$.  Let $C$ be  a right $H$-module and a coalgebra quotient of $C$ via a surjection $\pi: H\to C$.  View $A$ as a right $C$-comodule via the induced coaction $\varrho^A = (A\otimes_k \pi)\circ \bar{\varrho}^A$. Suppose that $C$ is a coseparable coalgebra (or that $k$ is algebraically closed and $C$ is a cosemisimple coalgebra). If  the (lifted) canonical map
$\tcan: A\ot_kA{\to} A\ot_k C$, $a\ot a'\mapsto a\varrho^A(a')$
is surjective, then $A$ is a principal $C$-extension of the coinvariants 
$S=A^{co C}$.
\label{cor.schsch}
\end{corollary}
\begin{proof}
This is an example of a Doi-Koppinen entwining structure $(A,C)_\psi$ with $\psi:c\ot a\mapsto \sum a\sw 0\ot ca\sw 1$, where $\bar{\varrho}^A(a) = \sum a\sw 0 \ot a\sw 1\in A\ot_k H$. Since $1_H$ is a group-like element in $H$, $e = \pi(1_H)$ is a group-like element in $C$. Note that
$$
\psi(e\ot a) = \sum a\sw 0\ot \pi(1_H)a\sw 1 = \sum a\sw 0\ot \pi(a\sw 1) = \roA(a),
$$
for $\pi$ is a right $H$-module map. Thus the $C$-coaction on $A$ has the required form. The map $\psi$ is bijective with the
inverse  $\psi^{-1} : A\ot_k C\to C\ot_k A$, 
$a\ot c\mapsto \sum c\,\mathsf{S}^{-1}(a\sw 1)\ot a\sw 0$, where $\mathsf{S}$ is the antipode in $H$. Any cosemisimple coalgebra over an algebraically closed field is coseparable, thus in either case all the assumptions of Theorem~\ref{thm.main} are satisfied and the assertion follows.
\end{proof}

An explicit, very important and geometrically interesting example of Galois corings (principal extensions) of the type described in Corollary~\ref{cor.schsch} has been recently contructed in \cite{BonCic:bij}. In this example $H$ is the Hopf algebra of functions on the quantum group $U_q(4)$, the algebra $A$ is the algebra of functions on the quantum 7-sphere $S_q^7$, the induced coalgebra $C$ is the coalgebra of functions on the quantum group $SU_q(2)$. Finally, the coinvariant algebra $S$ is the algebra of functions on the ``Etruscan" quantum 4-sphere $\Sigma_q^4$ introduced in \cite{BonCic:ins}.

 \section{Morphisms of corings and induced Galois and principal comodules} 
\label{sec.induce}

It frequently happens that there is a pair of corings related by a coring morphism, and one can prove that a comodule of one of these corings is, respectively, a Galois or a principal comodule. A question then arises, whether the induced comodule is also a Galois or a principal comodule. This is the main subject of the present section.

  \begin{lemma} 
  Let $(\gamma:\alpha) : (\cC:A)\to (\cD:B)$ be a morphism of corings. 
  Suppose that $M$ is a Galois right $\cC$-comodule and let $S = \Rend \cC 
  M$. If $M$ is a flat left $S$-module then for any $N\in\M^\cD$, the map 
  $$ 
  \vartheta_{M,N} : \Rhom \cD {M\otimes_AB} N \otimes_S M\to N\square_\cD 
  (B\otimes_A\cC), \;\; f\otimes m\mapsto \sum f(m\sw 0\otimes 
  1_B)\otimes m\sw 1 
  $$ 
  is an isomorphism of $k$-modules. This isomorphism is natural in $N$, 
  i.e., it is an isomorphism of functors $\Rhom \cD {M\ot_AB} - \ot_S M\simeq -\Box_\cC(B\ot_A\cC)$. 
  \label{lemma.theta} 
  \end{lemma} 
  \begin{proof} 
  Since $M$ is a Galois comodule and ${}_SM$ is flat, Theorem~\ref{thm.galoisc}(1) implies that $\cC$ is a flat left $A$-module. Thus 
  $-\square_\cD(B\otimes_A\cC)$ is a functor $\M^\cD\to \M^\cC$ and, 
  consequently, $N\square_\cD(B\otimes_A\cC)$ is a right 
  $\cC$-comodule. Therefore Theorem~\ref{thm.galoisc}(1) again implies that the evaluation map
  $$
\varphi_{N\square_\cD (B\otimes_A\cC)}: \Rhom \cC M {N\square_\cD (B\otimes_A\cC)}\ot_SM\to 
  N\square_\cD (B\otimes_A\cC), \quad  f\otimes m \mapsto f(m)
$$
 is an isomorphism of right $\cC$-comodules. 
  Combining $\varphi_{N\square_\cD (B\otimes_A\cC)}$ with   the general hom-tensor 
  relation  isomorphism of $k$-modules $\Rhom \cC M {N\square_\cD 
  (B\otimes_A\cC)} \simeq \Rhom \cD {M\otimes_A B} N$, we obtain the 
  desired isomorphism $ 
  \vartheta_{M,N}$. 

  To verify the naturality of maps $\vartheta_{M,N}$ we need to take any 
  morphism $g: N\to N'$ in $\M^\cD$ and establish that the following 
  diagram 
  $$\xymatrix{ 
      \Rhom \cD{M\ot_AB} N \ot_S  M \ar[rr] 
  ^(.6){\vartheta_{M,N}}\ar[d]_{ \Rhom \cD{M\ot_AB} g \ot_S  M}  & & 
          N\square_\cD(B\otimes_A\cC) \ar[d]^{g\square_\cD(B\ot_A\cC)} \\ 
       \Rhom \cD{M\ot_AB}{N'} \ot_S  M \ar[rr] ^(.6){\vartheta_{M,N'}}& 
  &N'\square_\cD(B\otimes_A\cC) 
  } $$ 
  is commutative. This follows immediately from the definition of 
  $\vartheta_{M,N}$. 
  \end{proof} 

  \begin{theorem} 
  Let $(\gamma:\alpha) : (\cC:A)\to (\cD:B)$ be a morphism of corings. 
  Suppose that $M$ is a Galois right $\cC$-comodule and let $S = \Rend \cC 
  M$ and $T=\Rend \cD {M\ot_A B}$.  If $M$ is a faithfully flat left 
  $S$-module then the following statements are equivalent. 
  \begin{zlist} 
  \item $\cD$ is a flat left $B$-module and $B\otimes_A\cC$ is a 
  faithfully coflat left $\cD$-comodule. 
  \item $M\ot_AB$ is a Galois right $\cD$-comodule and $M\ot_AB$  is a 
  faithfully flat left $T$-module. 
  \item $\cD$ is a flat left $B$-module and $\Rhom\cD{M\ot_AB}{-} 
  :\M^\cD\to \M_T$ is an equivalence with the inverse 
  $-\otimes_T(M\otimes_AB):\M_T\to\M^\cD$. 
  \item $\cD$ is a flat left $B$-module and $M\otimes_AB$ is a projective 
  generator in $\M^\cD$. 
  \end{zlist} 
  Furthermore, if $B$ is a quasi-Frobenius (QF) ring then the above 
  statements are equivalent to 
  \begin{zlist}\addtocounter{zlist}{4} 
  \item $\cD$ is a flat left $B$-module and $B\otimes_A\cC$ is an 
  injective cogenerator in ${}^\cD\M$. 
  \end{zlist} 
  \label{thm.Galois.main.1} 
  \end{theorem} 
  \begin{proof} 
  Clearly, the equivalences (2)$ \Lra$(3)$\Lra$(4) are contained in the 
  Galois comodule structure theorem. The equivalence of (1) and (5) (in 
  the case of a QF ring $B$) follows from the description of faithfully 
  coflat comodules for corings over QF rings (cf.\ 
  \cite[21.9(2)(ii)]{BrzWis:cor}). Thus it suffices to show that the statement 
  (1) is equivalent to the statement (3). 

  (3)$\Ra$(1) Since the functor $\Rhom\cD{M\ot_AB}{-}$ is an equivalence of Abelian categories,
  it reflects and preserves exact sequences. Furthermore, since ${}_SM$ is 
  a faithfully flat module, the composite functor $\Rhom\cD{M\ot_AB} 
  {-}\ot_SM$ preserves and reflects exact sequences. By 
  Lemma~\ref{lemma.theta} the functor $\Rhom\cD{M\ot_AB}{-}\ot_SM$ is 
  naturally isomorphic to the functor $-\square_\cD (B\ot_A\cC)$, thus the 
  latter also preserves and reflects exact sequences. Hence $B\ot_A\cC$ 
  is a faithfully coflat left $\cD$-comodule. 

  (1)$\Ra$(3) Take any right $T$-module $W$ and any right $\cD$-comodule 
  $N$. Since the functor $\Rhom\cD{M\ot_AB}{-}$ is right adjoint to 
  the functor $-\otimes_T (M\otimes_A B)$, there is an isomorphism 
  $$ 
  \rhom T W {\Rhom\cD{M\ot_AB}{N}} \to \Rhom \cD {W\ot_TM\ot_AB} N, 
  \quad f\mapsto F_f. 
  $$ 
  Explicitly, $F_f (w\ot m\ot b) = f(w)(m\ot b)$. We need to show that $f$ 
  is an isomorphism if and only if $F_f$ is an isomorphism. 

  First note that the map $\vartheta_{M,M\ot_AB}: T\otimes_SM=\Rend \cD 
  {M\otimes_AB}\otimes_SM \to (M\otimes_A B)\square_\cD(B\otimes_A\cC)$ in 
  Lemma~\ref{lemma.theta} is an isomorphism of left $T$-modules. Indeed, 
  take any $t, f\in T$ and $m\in M$ and compute 
  \begin{eqnarray*}
  \vartheta_{M,M\otimes_A B}(tf\ot m) &=&\sum tf(m\sw 0\ot 1_B)\ot m\sw 1 \\
&= &
  \sum t(f(m\sw 0\ot 1_B))\ot m\sw 1 = t\vartheta_{M,M\otimes_A B}(f\ot m), 
\end{eqnarray*}
  as required. Therefore we can construct  an isomorphism 
$$\theta: 
  W\ot_SM\to (W\ot_T(M\otimes_A B))\square_\cD(B\otimes_A\cC)$$ 
as a 
  composition 
  $$\xymatrix{ 
  W\ot_SM\ar[r]^{\simeq}& 
  M\ot_TT\ot_SM\ar[r]^(.4){W\ot_T\vartheta_{M,M\otimes_A B}} 
  &W\ot_T((M\otimes_A B)\square_\cD(B\otimes_A\cC)) &\\ 
    ~~~\ar[r]^(.2){\simeq}&(W\ot_T(M\otimes_A 
  B))\square_\cD(B\otimes_A\cC)& & &. 
  } 
  $$ 
  The last isomorphism is the consequence of the fact that $B\otimes_A\cC$ 
  is a (faithfully) coflat left $\cD$-comodule. In this way we are led to 
  the following commutative diagram 
  $$\xymatrix{ 
      W \ot_S  M \ar[rr] ^(.4){\theta}\ar[d]_{f \ot_S  M}  & & 
          (W\ot_T(M\otimes_A B))\square_\cD(B\otimes_A\cC) 
  \ar[d]^{F_f\square_\cD(B\ot_A\cC)} \\ 
       \Rhom \cD{M\ot_AB}N \ot_S  M \ar[rr] ^(.5){\vartheta_{M,N}}& 
  &N\square_\cD(B\otimes_A\cC) . 
  } $$ 
  Since the rows are ismorphisms, $M$ is a faithfully flat left $S$-module 
  and $B\otimes_A\cC$ is a faithfully coflat left $\cD$-comodule, the map 
  $f$ is an isomorphism if and only if $F_f$ is an isomorphism. Thus 
  $\Rhom \cD {M\ot_AB} -: \M^\cD\to M_T$ is an equivalence as required. 
  \end{proof} 

  The following lemma is an immediate consequence of the definition of a 
  faithfully coflat comodule, and gives one a criterion of the faithful 
  coflatness. 
  \begin{lemma} 
  Let $(\gamma:\alpha) : (\cC:A)\to (\cD:B)$ be a morphism of corings. 
  Suppose that: 
  \begin{blist} 
  \item $\cD$ is flat as a left $B$-module; 
  \item $B\ot_A\cC$ is a coflat left $\cD$-comodule; 
  \item the induced map $\tilde{\gamma}:B\ot_A\cC\to \cD$, $b\ot c\mapsto 
  b\gamma(c)$ is a split epimorphism in $\M^\cD$. 
  \end{blist} 
  Then $B\ot_A\cC$ is a faithfully coflat left $\cD$-comodule. 
  \label{lemma.criterion} 
  \end{lemma} 
  \begin{proof} 
  Condition (c) implies that $N\square_\cD (B\ot_A\cC) \neq 0$ for all 
  nonzero $N\in \M^\cD$, and thus the assertion follows from 
  \cite[21.7]{BrzWis:cor}. 
  \end{proof} 

The criterion of faithful flatness in Lemma~\ref{lemma.criterion}  assures that the principality of a comodule is carried over to the induced comodule.
  \begin{theorem} 
  Let $(\gamma:\alpha) : (\cC:A)\to (\cD:B)$ be a morphism of corings such 
  that the induced map $\tilde{\gamma}:B\ot_A\cC\to \cD$, $b\ot c\mapsto 
  b\gamma(c)$ is a split epimorphism in $\M^\cD$. Suppose that $M$ is a 
  principal $\cC$-comodule that is faithfully flat as a $k$-module,  $\cD$ is flat as a left $B$-module and 
  $B\ot_A\cC$ is a coflat left $\cD$-comodule. 
  Then $M\ot_AB$ is a principal $\cD$-comodule. 
  \label{theorem.projective} 
  \end{theorem} 
  \begin{proof} 
  Let $S=\Rend \cC M$. Since $\tilde{\gamma}$ is a split epimorphism in $\M^\cD$, there exists 
  a right $\cD$-comodule $V$ such that 
  $$ 
  \cD\oplus V\simeq B\ot_A\cC. 
  $$ 
  Since the cotensor product commutes with the colimits (cf.\ 
  \cite[21.3(3)]{BrzWis:cor}), the above isomorphism induces the 
  isomorphism of left $T$-modules 
  $$ 
  (M\ot_AB)\square_\cD\cD\oplus (M\ot_AB)\square_\cD V\simeq 
  (M\ot_AB)\square_\cD(B\ot_A\cC). 
  $$ 
  Note that $(M\ot_AB)\square_\cD\cD\simeq M\ot_A B$ and 
  $$ 
  (M\ot_AB)\square_\cD(B\ot_A\cC) \simeq \Rhom \cD {M\ot_A B}{M\ot_A 
  B}\ot_S M =T\ot_SM, 
  $$ 
  by Lemma~\ref{lemma.theta}. This is again an isomorphism of left 
  $T$-modules (compare the proof of Theorem~\ref{thm.Galois.main.1} 
  (1)$\Ra$(3)). Since $M$ is a projective left $S$-module, the induced 
  module $T\ot_S M$ is a projective left $T$-module. Thus ${}_T(M\ot_AB)$ 
  is a direct summand of a projective module, hence a projective module. 

  In view of Lemma~\ref{lemma.criterion} and Theorem~\ref{thm.fflat}, the hypothesis (1) in 
  Theorem~\ref{thm.Galois.main.1} is satisfied, hence $M\ot_A B$ is a Galois $\cD$-comodule. Since it is also a projective $T$-module, it is a principal comodule as 
claimed. 
  \end{proof} 

The $q$-deformed second Hopf fibration of \cite{BonCic:bij}, \cite{BonCic:ins} as described at the end of Section~\ref{sec.principal} is an example of induction of principal comodules provided by Theorem~\ref{theorem.projective}.

  \section{Duality and associated modules (non-commutative vector bundles)} 

The following theorem describes a remarkable duality of Galois comodules.
\begin{theorem}\label{thm.dual}~\\
\vspace{-\baselineskip}
\begin{zlist}
\item Let $M$ be a Galois right $\cC$-comodule and let $S=\Rend \cC M$. If $M$ is flat as a left $S$-module, then for any right $\cC$-comodule $W$,
$$
\Rhom \cC WM \simeq \rhom S {\Rhom \cC M W} S = (\Rhom \cC M W)^*,
$$
as left $S$-modules.
\item Let $N$ be a Galois left $\cC$-comodule and let $T=\Lend \cC N$. If $N$ is flat as a right $T$-module, then for any left $\cC$-comodule $V$,
$$
\Lhom \cC VN \simeq \lhom T {\Lhom \cC N V} T = {}^*(\Lhom \cC NV),
$$
as right $T$-modules.
\end{zlist}
\end{theorem}
\begin{proof} We only prove assertion (1), since (2) will follow by the right-left symmetry. Since ${}_SM$ is flat, the first part of the Galois comodule structure theorem, Theorem~\ref{thm.galoisc}, implies that $W\simeq \Rhom \cC M W\ot_S M$. Apply $\Rhom \cC - M$ to this isomorphism to deduce that
$$
\Rhom\cC W M \simeq \Rhom\cC{\Rhom \cC M W\ot_S M} M.
$$
Now, the hom-tensor relations \cite[18.10(2)]{BrzWis:cor} imply that
$$
\Rhom\cC{\Rhom \cC M W\ot_S M} M \simeq \rhom S{\Rhom \cC M W} {\Rhom \cC MM},
$$
i.e.,
$$
\Rhom \cC WM \simeq \rhom S {\Rhom \cC M W} S = (\Rhom \cC M W)^*,
$$
as required. Note that all the maps in this chain of isomorphism are maps of left $S$-modules.
\end{proof}

Since a dual of a Galois right $\cC$-comodule is a Galois left $\cC$-comodule and the endomorphism rings of these comodules are mutually isomorphic,  Theorem~\ref{thm.dual} leads immediately to the following
\begin{corollary}\label{cor1}
Let $M$ be a Galois right $\cC$-comodule and let $S=\Rend \cC M$. If $M^*$ is flat as a right $S$-module, then for any left $\cC$-comodule $V$,
$$
\Lhom \cC V{M^*} \simeq \lhom S {\Lhom \cC {M^*} V} S = {}^*(\Lhom \cC {M^*}V),
$$
as right $S$-modules.
\end{corollary}

Recall that for a right $\cC$-comodule $M$ that is finitely generated and projective as a right $A$-module, $\Rhom\cC M W \simeq W\Box_\cC M^*$, for any right $\cC$-comodule $W$ (cf.\ \cite[21.8]{BrzWis:cor}). Explicitly $f\mapsto \sum_i f(e^i)\ot \xi^i$, where $\{e^i\in M,\xi^i\in M^*\}_{i=1,\ldots, n}$ is a dual basis of $M$. Similar isomorphism holds for left $\cC$-comodules. Taking these isomorphisms into account, one obtains the following immediate consequence of Theorem~\ref{thm.dual} and Corollary~\ref{cor1}

\begin{corollary}\label{cor2}
With the notation and assumptions as in Theorem~\ref{thm.dual},
$$
\Rhom \cC WM \simeq  (W\Box_\cC M^*)^*, \qquad \Lhom \cC VN \simeq  {}^*({}^*N\Box_\cC V).
$$
Furthermore, for any Galois right $\cC$-comodule $M$ such that $M^*_S$ is flat, and for any left $\cC$-comodule $V$
$$
{}^*(M\Box_C V) \simeq \Lhom \cC V{M^*} .
$$
\end{corollary}

In particular, in the case of a Galois coring (i.e., when $A$ is a Galois comodule), $A^*\simeq {}^*A\simeq A$, and hence some of the stars can be removed in Corollary~\ref{cor2}, thus leading to
\begin{corollary}\label{cor3}
Let $\cC$ be a Galois $A$-coring, $g = \varrho^A(1)$ be the corresponding group-like element. Then  endomorphisms $S= \Rend \cC A$ come out as  $S = \{s\in A\; |\; sg=gs\}$ (cf.\ Example~\ref{ex.principal}). 

If $A$ is flat as a left $S$-module, then, for any right $\cC$-comodule $W$, there is an isomorphism of left $S$-modules
$$
(W_g^{co\;\cC})^* = \rhom S{W\Box_\cC A}S \simeq \Rhom \cC WA ,
$$
where $W_g^{co\;\cC} = \{ w\in W \; |\; \varrho^W(w) = w\ot g\}$ is a right $S$-module of  {\em $g$-coinvariants} of a right $\cC$-comodule $W$. Note that
$$\Rhom \cC WA = \{f\in \rhom AWA \; | \; \forall w\in W, \sum f(w\sw 0)w\sw 1 = gf(w)\}.
$$

If $A$ is flat as a right $S$-module, then, for any left $\cC$-comodule $V$, there is an isomorphism of right $S$-modules
$$
{}^*({}^{co\;\cC}V_g) = \lhom S{A\Box_\cC V}S \simeq \Lhom \cC VA ,
$$
where ${}^{co\;\cC} V_g= \{ v\in V \; |\; {}^V\varrho(v) = g\ot v\}$ is a left $S$-module of  {\em $g$-coinvariants} of a left $\cC$-comodule $V$. Note that
$$
\Lhom \cC VA = \{f\in \lhom AVA \; | \; \forall v\in V, \sum v\sw{-1} f(v\sw 0) = f(v)g\} .
$$
\end{corollary}
\begin{example}\label{ex.sections}
 As an example for Corollary~\ref{cor3}, take a coalgebra-Galois $C$-extension $S\subseteq A$ over a field $k$, i.e., $A$, $C$ and $S$ are as  in Example~\ref{ex.principal} but only  the condition Example~\ref{ex.principal} (1) is required to hold. As recalled in the proof of Example~\ref{ex.principal}, $\cC= A\ot_kC$ is then an $A$-coring, and $A$ is a Galois comodule (hence $A\ot_k C$ is a Galois coring). The group-like element is $\varrho^A(1) = \sum 1\sw 0\ot1\sw 1\in A\ot_k C$. Take a left $C$-comodule $U$. Then $V=A\ot_kU$ is a left $\cC$-comodule with the coaction $a\ot u\mapsto \sum a\ot u\sw{-1}\ot u\sw 0\in A\ot_kC\ot_kU\simeq A\ot_kC\ot_AA\ot_kU$, and $\Lhom \cC V A = \hom\psi U A$, where
$$
\hom\psi U A  = \{f\!\in \!\hom k U A \,|\, \forall u\in U, \sum \psi(u\sw{-1}\ot f(u\sw 0)) = \sum f(u)1\sw 0\ot 1\sw 1\}.
$$
Here $\psi$ is the canonical entwining map (cf.\  Example~\ref{ex.principal}(3)). On the other hand
${}^{co\;\cC}V_g = A\Box_C U$. Thus, if $A_S$ is flat, Corollary~\ref{cor3} implies that $\Lhom S {A\Box_C U} S \simeq \hom\psi U A$ as right $S$-modules and we obtain (a part of) \cite[Theorem~4.3]{Brz:mod}.

Now take a right $C$-comodule $X$. Then $W=X\ot_k A$ is a right $\cC$-comodule with the coaction
$\varrho^W: x\ot a\mapsto \sum x\sw 0\ot \psi(x\sw 1\ot a)$. In this case
$$
\Rhom \cC W A \simeq \Rhom C XA
$$
and 
$$
 W_g^{co\;\cC} = (X\ot_k A)_0 := \{ x\ot a\in X\ot_k A \; | \; \sum x\sw 0\ot \psi(x\sw 1\ot a) = \sum x\ot a1\sw0\ot 1\sw 1\}.
$$
Then, if $A$ is a flat left $S$-module, Corollary~\ref{cor3} yields  the isomorphism of left $S$-modules $\rhom S {(X\ot_kA)_0}S\simeq \Rhom C X A$ in \cite[Theorem~5.4]{Brz:mod}.
\end{example}

As explained in \cite{Brz:mod}, both $A\Box_C U$ and $(X\ot_k A)_0$ in Example~\ref{ex.sections} have the non-commutative geometric meaning of fibre bundles associated to non-commutative principal bundles. The isomorphisms described in this example generalise the classical correspondence between sections of a fibre bundle and covariant (tensorial) functions on the principal bundle with values in the fibre (functions of type $\rho$). Thus Theorem~\ref{thm.dual} (and its corollaries) can be understood as an algebraic origin of this deep geometric fact.

One of the main motivations for introducing principal extensions in \cite{BrzHaj:che} is the observation that  if $A$ is a principal $C$-extension and $X$ is a finite dimensional vector space, then left $S$-module $\Rhom C X A$ is a finitely generated projective module, hence it can be truly interpreted as a module of sections on a non-commutative vector bundle in the sense of Connes (cf.\ \cite{Con:non}). In particular, in this way one can study some aspects of the $K$-theory of principal extensions (Chern-Galois characters). From the non-commutative geometry point of view it is therefore extremely interesting to study what additional properties must be imposed on a principal $\cC$-comodule $M$ to make $\Rhom \cC M W$ a finitely generated left $S$-module for any right $\cC$-comodule $W$ that is finitely generated and projective as an $A$-module. One of the possibilities is explored in the following
\begin{proposition}\label{prop.proj}
Let $k$ be a field, $M$ be a Galois right $\cC$-comodule  with the endomorphism ring  $S=\Rend\cC M$. 
View $S\ot_k M$ as a right $\cC$-comodule via $S\ot_k\roM$ and 
$M^*\ot_k S$ as a left $\cC$-comodule via ${}^{M^*}\!\varrho \ot_k S$.
\begin{zlist}
\item 
If ${}_SM$ is faithfully flat and there exists a right $S$-module left $\cC$-comodule section of the action  
$M^*\ot_k S\to M^*$, then, for every left $\cC$-comodule $V$ that is finitely generated 
projective as a left $A$-module, $\Lhom \cC V{M^*} $ is a finitely generated projective 
right $S$-module.
\item 
If $M^*_S$ is faithfully flat and there exists a left $S$-module right $\cC$-comodule section of the action  
$S\ot_k M\to M$, then, for every right $\cC$-comodule $W$ that is finitely generated 
projective as a right $A$-module, $\Rhom \cC W{M} $ is a finitely generated projective 
left $S$-module.
\end{zlist}
\end{proposition}
\begin{proof}
(1) 
If a left $\cC$-comodule $V$ is a 
finitely generated projective left $A$-module then ${}^*V$ is a right $\cC$-comodule 
and $\Lhom \cC V{M^*}\simeq {}^*V\Box_\cC M^*$ as right $S$-modules. Thus 
one can consider the following chain of isomorphisms of right $A$-modules
$$
\Lhom \cC V{M^*}\ot_SM\simeq ({}^*V\Box_\cC M^*)\ot_S M \simeq {}^*V\Box_\cC (M^*\ot_S M) \simeq {}^*V\Box_\cC\cC \simeq {}^*V.
$$
The second isomorphism is the consequence of the fact that ${}_SM$ is flat, the third one follows from the bijectivity of the 
canonical map $\can_M$. Since $V$ is a finitely generated left $A$-module, its dual is a finitely generated right $A$-module. Thus $\Lhom \cC V{M^*}\ot_SM$ is a finitely generated right $A$-module. Therefore,  there exists a finite number of elements $x_{ij}\in\Lhom \cC V{M^*}$ such that  $\Lhom \cC V{M^*}\ot_SM$ is generated by $\sum_j x_{ij}\ot e^j$, where  $e^j$ are generators  of $M_A$. Let $X$ be a submodule of $\Lhom \cC V{M^*}$ generated by the $x_{ij}$. Applying $-\otimes_S M$ to the inclusion $0\to X\to \Lhom \cC V{M^*}$ one obtains a right $A$-module surjection $X\ot_SM\to  \Lhom \cC V{M^*}\ot_SM\to 0$. Since ${}_SM$ is faithfully flat, also the inclusion $X\to \Lhom \cC V{M^*}$ is a surjection, and since $X$ is a finitely generated right $S$-module, so is  $\Lhom \cC V{M^*}$ (cf.\   proof of \cite[Ch.\ 1\S 3 Prop.\ 11]{Bou:com}). 

Let $\sigma: M^*\to M^*\ot_k S$ be a right $S$-linear left $\cC$-colinear section of the action $M^*\ot_k S\to M^*$. Write $\theta: \Lhom \cC V{M^*}\to {}^*V\Box_\cC M^*$ for the isomorphism of right $S$-modules  $\Lhom \cC V{M^*}\simeq {}^*V\Box_\cC M^*$. Then $\sigma_V = (\theta^{-1}\ot _kS)\circ ({}^*V\Box_\cC \sigma)\circ \theta$ is a right $S$-module section of the multiplication map $\Lhom \cC V{M^*}\ot_k S\to \Lhom \cC V{M^*}$. Thus $\Lhom \cC V{M^*}$ is a projective module. This completes the proof of the assertion that $\Lhom \cC V{M^*}$ is a finitely generated projective right $S$-module.

(2) Follows by similar arguments as part (1).
\end{proof}

Unfortunately, the conditions in Proposition~\ref{prop.proj} are too strong to cover the case of 
principal extensions. Thus the problem of finding suitable conditions remains open. Another possible criterion is given in Proposition~\ref{prop.split.proj}.

Theorem~\ref{thm.dual} implies also the following reflexivity property of Galois comodules.
\begin{corollary}\label{cor4}
Let $M$ and $N$ be Galois right $\cC$-comodules and set $S= \Rend \cC M$ and $T=\Rend \cC N$. If ${}_SM$ and ${}_TN$ are flat, then
$$
\Rhom \cC M N \simeq \rhom T {\rhom S {\Rhom \cC M N} S}T = ((\Rhom \cC M N)^*)^*,
$$
as $(T,S)$-bimodules.
\end{corollary}
\begin{proof}
Apply Theorem~\ref{thm.dual} twice, and notice that the isomorphism in Theorem~\ref{thm.dual},
$$
\Rhom \cC NM \to \rhom S{\Rhom \cC MN}S, \quad f\mapsto [\phi\mapsto f\circ\phi],
$$
is an $(S,T)$-bimodule map.
\end{proof}

\section{Relatively injective Galois comodules}
For any right $\cC$-comodule $M$, the comodule endomorphism ring $S=\Rend \cC M$ is a subring of the module endomorphism ring $\hat{S} = \rend AM$. Thus there is a ring extension $S\to \hat{S}$. The properties of this extension capture properties of the comodule $M$. This is most profound in the case of Galois comodules. 

Recall that a right $\cC$-comodule $M$ is called a  {\em $(\cC, A)$-injective comodule} or an {\em $A$-relatively injective comodule} if,  for every $\cC$-comodule map $i:N \to L$ that has
a  retraction in $\cM_A$, every diagram 
$$ 
\xymatrix{ N \ar[rr]^i \ar[rd]_f & & L \\ 
	      &   M   &   } 
$$ 
in $\cM^\cC$ can be completed commutatively by some $g:L\to M$ in $\cM^\cC$. Equivalently, $M$ is 
a   $(\cC, A)$-injective comodule if the coaction $\roM : M\to M\ot_A \cC$ has a right $\cC$-comodule retraction $\pi_M: M\ot_A\cC\to M$. Here $M\ot_A\cC$ is a right $\cC$-comodule with the coaction $M\ot_A\DC$. $(\cC,A)$-injective left comodules are defined in a symmetric way.

Recall also that a ring extension  $S\to R$ is called a {\em split extension} provided there exists an $(S,S)$-bimodule map $\sigma: R\to S$ such that $\sigma(1_R) =1_S$.  We begin with the following simple
\begin{lemma}\label{lemma.s}~\\
\vspace{-\baselineskip}
\begin{zlist}
\item Let $M$ be a Galois right $\cC$-comodule and let $S=\Rend \cC M$. Then
$$
\Rhom \cC \cC M \simeq \rhom S {M^*} S ,
$$
as left $S$-modules.
\item Let $N$ be a Galois left $\cC$-comodule and let $T=\Lend \cC N$. Then 
$$
\Lhom \cC \cC N \simeq \lhom T {N ^*} T ,
$$
as right $T$-modules.
\end{zlist}
In particular, for a a Galois right $\cC$-comodule $M$, $\Lhom \cC \cC {M^*} \simeq \lhom S {M} S$.
\end{lemma}
\begin{proof}
We only prove (1) as (2) can be obtained by the left-right symmetry. This is a consequence of the following chain of isomorphisms:
\begin{eqnarray*}
\Rhom \cC \cC M &\simeq& \Rhom \cC {M^*\ot_S M} M \\
& \simeq& \rhom S {M^*} {\Rhom \cC MM} \simeq \rhom S {M^*} S .
\end{eqnarray*}
The first isomorphism is obtained by applying $\Rhom \cC - M$ to the canonical map $\can_M$, while the second is the hom-tensor relation for modules and comodules (cf.\ \cite[18.10(2)]{BrzWis:cor}). The final assertion follows from the observation that $M^*$ is a Galois left $\cC$-comodule and from assertion (2).
\end{proof}
 
\begin{theorem}\label{thm.split} ~\\
\vspace{-\baselineskip}
\begin{zlist}
\item Let $M$ be a Galois right $\cC$-comodule,  let $S=\Rend \cC M$ and $\hat{S} = \rend AM$. Then $M$ is a $(\cC,A)$-injective comodule if and only if there exists a right  $S$-module map $\sigma:\hat{S}\to S$ such that $\sigma(1_{\hat{S}}) = 1_S$.

\item Let $N$ be a Galois left $\cC$-comodule, let $T=\Lend \cC N$ and $\hat{T} = \lend AN$. Then $N$ is a $(\cC,A)$-injective comodule if and only if there exists a left  $T$-module map $\tau:\hat{T}\to T$ such that $\tau(1_{\hat{T}}) = 1_T$.

\item Let $M$ be a Galois right $\cC$-comodule, let $S=\Rend \cC M$ and $\hat{S} = \rend AM$. The following statements are equivalent
\begin{blist}
\item  $S\to \hat{S}$ is a split extension;
\item the right coaction $\roM$ has a right $\cC$-comodule left $S$-module retraction;
\item the left coaction ${}^{M^*\!\!}\varrho$ has a left $\cC$-comodule right $S$-module retraction.
\end{blist}
In particular, if $S\to \hat{S}$ is a split extension, then $M$ is a $(\cC,A)$-injective right $\cC$-comodule and $M^*$ is a $(\cC,A)$-injective left $\cC$-comodule.
\end{zlist}
\end{theorem}
\begin{proof}
(1) Consider  the following chain of isomorphisms
\begin{eqnarray*}
\rhom S {\hat{S}} S &\simeq & \rhom S {M\ot_A M^*} S \simeq \rhom A M {\rhom S {M^*} S}\\
&\simeq&\rhom A M {\Rhom \cC \cC M} \simeq \Rhom \cC {M\ot_A\cC} M. 
\end{eqnarray*}
The first isomorphism is the consequence of the canonical isomorphism $\rend AM \simeq M\ot_A M^*$ that holds for any finitely generated projective module. The second and the fourth isomorphisms are the hom-tensor relations for modules and comodules respectively. The third isomorphism is obtained by applying $\rhom A M -$ to the isomorphism in Lemma~\ref{lemma.s}. Explicitly, the composite isomorphism $\Theta: \rhom S {\hat{S}} S \to  \Rhom \cC {M\ot_A\cC} M$ comes out as follows. First, for any $c\in \cC$ write
$$
\can_M^{-1}(c) = \sum c\su 1\ot c\su 2\in M^*\ot_S M.
$$
Then, for all $\sigma \in \rhom S {\hat{S}} S$, $m\in M$, $c\in \cC$, $\pi\in \Rhom \cC {M\ot_A\cC} M$ and $\hat{s}\in\hat{S}$, the isomorphism $\Theta$ and its inverse $\Theta^{-1}$ read
$$
\Theta(\sigma)(m\ot c) = \sum \sigma(mc\su 1(-))(c\su 2), \qquad \Theta^{-1}(\pi)(\hat{s}) = \pi\circ(\hat{s}\ot_A\cC)\circ\roM.
$$
The fact that $\Theta$ and $\Theta^{-1}$ are mutual inverses can also be shown directly as follows. First note that $\can_M^{-1}(c) $ has the following properties:
$$
\sum c\su 1(c\su 2\sw 0)c\su 2\sw 1 = c, \qquad \forall c\in \cC, \eqno{(A)}
$$
for $(\can_M\circ\can_M^{-1})(c) = c$, and
$$
\xi\ot_S m = \sum \xi(m\sw 0)m\sw 1\su 1\ot_S m\sw 1\su 2 , \qquad \forall m\in M,\; \xi\in M^*, \eqno{(B)}
$$
for $(\can_M^{-1}\circ\can_M)(\xi\ot m) = \xi\ot m$. Let $\{e^i\in M,\xi^i\in M^*\}_{i=1,\ldots, n}$ be the dual basis of $M_A$. In view of the isomorphism $\hat{S}\simeq M\ot_A M^*$, $\hat{s}\mapsto \sum_i \hat{s}(e^i)\ot \xi^i$, $m\ot\xi\mapsto m\xi( -)$, the property (B) implies that, for all $\hat{s}\in\hat{S}$ and $m\in M$,
$$
\hat{s}\ot_S m = \sum \hat{s}(m\sw 0)m\sw 1\su 1(-)\ot_S m\sw 1\su 2. \eqno{(C)}
$$
With these equalities at hand  we can compute, for all $\sigma \in \rhom S{\hat{S}}S$, $\hat{s}\in\hat{S}$ and $m\in M$,
\begin{eqnarray*}
(\Theta^{-1}\circ\Theta)(\sigma) (\hat{s})(m) &=&\Theta(\sigma)( \sum \hat{s}(m\sw 0)\ot m\sw 1)\\
& = &  \sum\sigma( \hat{s}(m\sw 0)m\sw 1\su 1(-))( m\sw 1\su 2) = \sigma(\hat{s})(m),
\end{eqnarray*}
by the property (C). On the other hand the use of property (A) entails that, for all $\pi\in \Rhom \cC {M\ot_A\cC} M$, $m\in M$ and $c\in \cC$,
\begin{eqnarray*}
(\Theta\circ\Theta^{-1})(\pi)(m\ot c) &=& \sum \Theta^{-1}(\pi)(mc\su 1(-))(c\su 2)\\
&=& \sum\pi(m\ot c\su 1(c\su 2\sw 0)c\su 2\sw 1) = \pi(m\ot c).
\end{eqnarray*}
Thus $\Theta$ and $\Theta^{-1}$ are mutual inverses as claimed on the basis of the chain of isomorphisms displayed at the beginning of the proof.

Suppose that $\sigma \in \rhom S{\hat{S}}S$ has the property $\sigma(1_{\hat{S}}) =1_S$. Then, for all $m\in M$,
$$
\Theta(\sigma)(\sum m\sw 0\ot m\sw 1) = \sum \sigma(m\sw 0m\sw 1\su 1(-))( m\sw 1\su 2) = \sigma(1_{\hat{S}})(m) = m,
$$
where the second equality follows from  property (C), by setting $\hat{s} = 1_{\hat{S}}$. Thus $\Theta(\sigma)$ is a retraction of $\roM$, and hence $M$ is a $(\cC,A)$-injective comodule. Conversely, if $M$ is a $(\cC,A)$-injective comodule and $\pi_M$ is a retraction of $\roM$, then, for all $m\in M$,
$$
\Theta^{-1}(\pi_M)(1_{\hat{S}})(m) =\pi_M( \sum  m\sw 0\ot m\sw 1) = m,
$$
i.e., $\Theta^{-1}(\pi_M)(1_{\hat{S}}) = 1_S$, as required. 

The assertion (2) is proven in the similar way to the proof of (1).

(3) The equivalence (a) $\Lra$ (b) follows from the observation that the maps $\Theta$ and $\Theta^{-1}$, constructed in the proof of assertion (1), preserve the left $S$-linearity. Indeed, if $\sigma \in \rhom S {\hat{S}} S$ is left $S$-linear then for all $s\in S$, $m\in M$ and $c\in \cC$,
\begin{eqnarray*}
\Theta (\sigma)(s(m)\ot c) &=& \sum \sigma(s(m) c\su 1(-))(c\su 2) = \sum \sigma(s(m c\su 1(-)))(c\su 2)\\ 
&=& \sum s \sigma(m c\su 1(-))(c\su 2) = s\Theta (\sigma)(m\ot c) ,
\end{eqnarray*}
so that $\Theta(\sigma)$ is also left $S$-linear. On the other hand, $\pi\in \Rhom \cC {M\ot_A\cC} M$ is left $S$-linear provided for all $s\in S$, $\pi\circ (s\ot_A\cC) = s\circ \pi$. If this is so, then the definition of $\Theta^{-1}(\pi)$ immediately implies that $\Theta^{-1}(\pi) $ is also left $S$-linear. In view of this and part (1), the coaction has a right $\cC$-colinear left $S$-linear retraction if and only if there exists $\sigma \in \lrhom SS{\hat{S}}S$ such that $\sigma(1_{\hat{S}}) = 1_S$, i.e., if and only if $S\to \hat{S}$ is a split extension.

The equivalence (a) $\Lra$ (c) is proven in a similar manner by noting that $M^*$ is a Galois left $\cC$-comodule, and $\lend A{M^*}\simeq \rend AM$ and $\Lend \cC {M^*} \simeq \Rend \cC M$. The final assertion is obvious. 
\end{proof}

Theorem~\ref{thm.split} leads to the following conditions for a Galois right comodule $M$ to be  faithfully flat as a left $S$-module. 

\begin{proposition}\label{prop.ff}
Let $M$ be a Galois right $\cC$-comodule and let $S=\Rend \cC M$ and $\hat{S} = \rend AM$. If either
\begin{blist}
\item ${}_SM$ is  flat  and there exists $\sigma\in\lhom S{\hat{S}}S$ such that $\sigma(1_{\hat{S}}) = 1_S$, or
\item ${}_A\cC$ is  flat  and $S\to \hat{S}$ is a split extension,
\end{blist}
then $M$ is a faithfully flat left $S$-module.
\end{proposition}
\begin{proof} (a) If ${}_SM$ is flat, then ${}_A\cC$ is flat and the counit of adjunction (the evaluation map) $\varphi_N : \Rhom \cC MN\ot_SM\to N$ is bijective for any right $\cC$-comodule $N$ by Theorem~\ref{thm.galoisc}(1). Since $M$ is a Galois right $\cC$-comodule, the dual module $M^*$ is a Galois left $\cC$-comodule. Furthermore, $\lend A{M^*}\simeq \rend AM$ and $\Lend \cC {M^*} \simeq \Rend \cC M$ via the isomorphism $\Gamma_M$ described in equation (\ref{Gamma}), and thus the existence of $\sigma$ implies that $M^*$ is a $(\cC,A)$-injective comodule by Theorem~\ref{thm.split}(2). Let $\pi_{M^*}\in \Lhom \cC {\cC\ot_A M^*} {M^*}$ be a retraction of the left coaction ${}^{M^*\!\!}\varrho$ corresponding to $\sigma$ as in Theorem~\ref{thm.split}. Taking into account $\Gamma_M$ (so that the elements of $S$ and $\hat{S}$ are understood as maps on $M$) and the relationship between $\cC$-coactions on $M$ and $M^*$ described in equation (\ref{left.coa1}), this correspondence has the following explicit form, for all $\hat{s}\in\hat{S}$, 
$$
\sigma(\hat{s}) = \sum_i \hat{s}(e^i)\sw 0\pi_{M^*}(\hat{s}(e^i)\sw 1\ot \xi^i)(-), \eqno{(*)}
$$
where $e^i\in M$, $\xi^i \in M^*$ is the finite dual basis of $M$. Let, for all $X\in \M_S$,
$$
\nu_X : X\to \Rhom \cC M {X\ot_SM}, \qquad x\mapsto [m\mapsto x\ot m]
$$
denote the unit of the adjunction. For any $f\in \Rhom \cC M {X\ot_SM}$ and $m\in M$, write $f(m) = \sum f(m)\suc 1\ot f(m)\suc 2\in X\ot_S M$. We claim that the map
$$
\nu_X^{-1} : \Rhom \cC M {X\ot_SM}\to X, \qquad f\mapsto \sum_i f(e^i)\suc 1\sigma(f(e^i)\suc 2\xi^i(-)),
$$
is the inverse of $\nu_X$. Indeed, for all $x\in X$,
\begin{eqnarray*}
(\nu_X^{-1}\circ\nu_X)(x) &=& \sum_i \nu_X(x)(e^i)\suc 1\sigma (\nu_X(x)(e^i)\suc 2\xi^i(-)) \\
&= &\sum_ix\sigma( e^i\xi^i(-))=x\sigma(1_{\hat{S}}) =x.
\end{eqnarray*}
On the other hand, in the view of the correspondence ($*$) and the fact that $f$ is a right $\cC$-comodule map,
\begin{eqnarray*}
\nu_X^{-1}(f) &=& \sum_i f(e^i)\suc 1(f(e^i)\suc 2\sw 0\pi_{M^*}(f(e^i)\suc 2\sw 1\ot \xi^i)(-)) \\
&=& \sum_i f(e^i\sw 0)\suc 1(f(e^i\sw 0)\suc 2\pi_{M^*}(e^i\sw 1\ot \xi^i)(-)) .
\end{eqnarray*}
Therefore, for all $m\in M$, $f\in \Rhom \cC M {X\ot_SM}$,
\begin{eqnarray*}
(\nu_X\circ\nu_X^{-1})(f)(m) &=& \sum_i f(e^i\sw 0)\suc 1(f(e^i\sw 0)\suc 2\pi_{M^*}(e^i\sw 1\ot \xi^i)(-))\ot_S m\\
&=& \sum_i f(e^i\sw 0)\suc 1\ot _S f(e^i\sw 0)\suc 2\pi_{M^*}(e^i\sw 1\ot \xi^i)(m)\\
&=& \sum_i f(e^i)\suc 1\ot _S f(e^i)\suc 2\pi_{M^*}(\xi^i\sw {-1}\ot \xi^i\sw 0)(m) \qquad \mbox{(by eq.\ (\ref{left.coa.2}))}\\
&=& \sum_i f(e^i)\xi^i(m) = f(m),
\end{eqnarray*}
where the penultimate equality follows from the fact that $\pi_{M^*}$ is a retraction of the coaction of $M^*$. 

In this way we have proven that $\Rhom \cC M - :\M^\cC\to \M_S$ is an equivalence. We have already observed that $\cC$ is a flat left $A$-module. Now Theorem~\ref{thm.galoisc}(2) implies that $M$ is a faithfully flat left $S$-module.

(b) Since $S\to \hat{S}$ is a split extension,  $M^*$ is a $(\cC,A)$-injective left $\cC$-comodule by Theorem~\ref{thm.split}. Furthermore, the retraction of the left coaction ${}^{M^*\!\!}\varrho$ is a right $S$-module map. Hence, for all $N\in \M^\cC$,  there is an isomorphism (the tensor-cotensor relation)
$$
N\Box_{\cC} (M^*\ot_S M) \simeq (N\Box_\cC M^*)\ot_S M,
$$
(cf.\ \cite[21.4, 21.5]{BrzWis:cor}). Furthermore, $M_A$ is finitely generated projective, so $N\Box_\cC M^*\simeq \Rhom \cC MN$, and therefore $
N\Box_{\cC} (M^*\ot_S M) \simeq  \Rhom \cC MN\ot_S M$. Thus one can consider the following commutative diagram
$$
\xymatrix{
0 \ar[r]& \Rhom \cC M N \ot_S M \ar[r]\ar[d]^{\varphi_N} & N\ot_A M^*\ot_S M \ar[r]\ar[d]^{N\ot_A\can _M} & N\ot_A\cC \ot_A M^*\ot_S M \ar[d]^{N\ot_A\cC\ot_A\can _M}\\
0 \ar[r] & N\ar[r]^{\varrho^N} & N\ot_A\cC \ar[r] & N\ot_A\cC\ot_A\cC}
$$
The top row is the defining sequence of $N\Box_\cC (M^*\ot_S M)\simeq  \Rhom \cC MN\ot_S M$, hence it is an exact sequence. The last map in the bottom row  is $\varrho^N\ot_A\cC - N\ot_A \DC$, hence the sequence is exact by the coassociativity of the coaction. Since $\can_M$ is an isomorphism, so are $N\ot_A\can_M$ and $N\ot_A\cC\ot_A\can_M$. Thus $\varphi_N$ is an isomorphism for any right $\cC$-comodule $N$. By the same arguments as in part (a) we conclude that $\Rhom\cC M - $ is an equivalence. Since ${}_A\cC$ is flat, Theorem~\ref{thm.galoisc}(2) implies that $M$ is a faithfully flat left $S$-module.
\end{proof}

A special case of Proposition~\ref{prop.ff} is obtained by taking $A$  a coalgebra-Galois $C$-extension, $\cC=A\ot_kC$ and $M=A$. In this case one derives \cite[Proposition~4.4]{Brz:mod}. Furthermore, if $A$ is a Hopf-Galois $H$-extension, $\cC=A\ot_kH$ and $M=A$, Proposition~\ref{prop.ff} gives \cite[Theorem~2.11]{DoiTak:Hop}. 

Combined with Theorem~\ref{thm.split}, the criterion of faithful flatness in Proposition~\ref{prop.ff}(b) leads also to the following criterion for modules associated to a Galois comodule to be finitely generated projective modules.
\begin{proposition}\label{prop.split.proj}
Let $k$ be a field, and let $M$ be a Galois right $\cC$-comodule, $S=\Rend \cC M$ and $\hat{S} = \rend AM$. Suppose that $S\to \hat{S}$ is a split extension. Then:
\begin{zlist}
\item 
If $M^*_S$ is projective (i.e., $M^*$ is a principal left comodule) and ${}_A\cC$ is flat, then, for any left $\cC$-comodule $V$ that is finitely generated 
projective as a left $A$-module, $\Lhom \cC V{M^*} $ is a finitely generated projective 
right $S$-module.
\item 
If ${}_SM$ is projective (i.e., $M$ is a principal right comodule) and $\cC_A$ is flat, then, for any right $\cC$-comodule $W$ that is finitely generated 
projective as a right $A$-module, $\Rhom \cC W{M} $ is a finitely generated projective 
left $S$-module.
\end{zlist}
\end{proposition}
\begin{proof}
(1) By Proposition~\ref{prop.ff}, $M$ is a faithfully flat left $S$-module. Hence by the same means as in the proof of Proposition~\ref{prop.proj} one proves  that $\Lhom \cC V{M^*} $ is a finitely generated right $S$-module.

Let $\sigma:M^*\to M^*\ot_k S$ be a right $S$-module section of the multiplication map. 
For all $\xi\in M^*$, write $\sigma(\xi) =\sum \xi\suc 1\ot \xi\suc 2\in M^*\ot_k S$, so that, for all $m\in M$, $\sum \xi\suc 1(\xi\suc 2(m)) =\xi(m)$. Let $e_i\in V$, $\eta_i\in {}^*V$ be a dual basis of $V$. Define a right $S$-linear map 
$\sigma_V:  \lhom AV {M^*}\to \lhom AV {M^*}\ot_k S$, by
$$
f \mapsto \sum_i \eta_i(-)\sigma(f(e_i))=  \sum_i \eta_i(-)f(e_i)\suc 1\ot f(e_i)\suc 2.
$$
Then, for all $v\in V$ and $m\in M$,
\begin{eqnarray*}
\sum_i ((\eta_i(-)f(e_i)\suc 1)\cdot f(e_i)\suc 2)(m)(v) &=& \sum_i \eta_i(v)f(e_i)\suc 1(f(e_i)\suc 2(m))\\
&=& \sum_i \eta_i(v)f(e_i)(m)\\
& = &\sum_i f(\eta_i(v)e_i)(m) = f(v)(m).
\end{eqnarray*}
This means that $\sigma_V$ is  a section of the product $\lhom A V{M^*}\ot_k S\to \lhom A V{M^*}$, and hence $\lhom A V{M^*}$ is a projective right $S$-module. 

By Theorem~\ref{thm.split}(3) there is a right $S$-linear, left $\cC$-colinear retraction $\pi: \cC\ot_A {}^*M\to {}^*M$ of the coaction ${}^{M^*}\!\varrho$. We claim that the map 
$$
\pi_V: \lhom A V{M^*}\to \Lhom \cC V{M^*}, \qquad f\mapsto  \pi \circ (\cC\ot_A f)\circ {}^V\!\varrho,
$$ 
is a right $S$-linear retraction of the defining inclusion $\Lhom \cC V{M^*}\subseteq \lhom A V{M^*}$. Indeed, $\pi_V(f)$ is left $\cC$-colinear, since it is a composition of left $\cC$-colinear maps. Furthermore, for all $f\in \lhom A V{M^*}$, $s\in S$ and $v\in V$,
\begin{eqnarray*}
\pi_V(fs)(v) &=& \sum\pi(v\sw{-1}\ot (fs)(v\sw 0)) = \sum\pi(v\sw{-1}\ot f(v\sw 0)\circ s) \\
&=&   \sum\pi(v\sw{-1}\ot f(v\sw 0))\circ s = (\pi_V(f)s)(v),
\end{eqnarray*}
where the penultimate equality is a consequence of the right $S$-linearity of $\pi$. Thus $\pi_V$ is a right $S$-linear map. Finally, if $f\in \Lhom \cC V{M^*}$, then 
$$
\pi \circ (\cC\ot_A f)\circ {}^V\!\varrho = \pi \circ{}^{M^*}\!\varrho\circ f =f,
$$ 
since $\pi$ is a section of the coaction. Therefore, $\Lhom \cC V{M^*}$ is a direct summand of a projective right $S$-module, hence a projective module.
This completes the proof that  $\Lhom \cC V{M^*}$ is a finitely generated projective right $S$-module, as required.

(2) This is proven in the analogues way to the proof of part (1). 
\end{proof}

  \section*{Acknowledgements} 
  I would like to thank Jos\'e G\'omez-Torrecillas and Bodo Pareigis for discussions. I would also like to thank the Engineering and Physical 
  Sciences Research 
  Council for an Advanced Fellowship.

  \end{document}